\documentstyle[12pt]{article}         
\newfont{\twlvmsb}{msbm10 scaled\magstep1}
\newfont{\ninemsb}{msbm9}
\newfont{\sixmsb}{msbm6}
\newfam\msbfam
\textfont\msbfam=\twlvmsb
\scriptfont\msbfam=\ninemsb
\scriptscriptfont\msbfam=\sixmsb
\catcode`\@=11
\def\Bbb{\ifmmode\let\next\Bbb@\else
  \def\next{\errmessage{Use \string\Bbb\space only in math mode}}\fi\next}
\def\Bbb@#1{{\Bbb@@{#1}}}
\def\Bbb@@#1{\fam\msbfam#1}
\newfont{\largeeufm}{eufm10 scaled\magstep4}
\newfont{\twlveufm}{eufm10 scaled\magstep1}
\newfont{\elveufm}{eufm10 at 11pt}
\newfont{\teneufm}{eufm10}
\newfont{\nineeufm}{eufm9}
\newfam\eufam
\textfont\eufam=\twlveufm
\scriptfont\eufam=\nineeufm
\def\frak{\ifmmode\let\next\frak@\else
\def\next{\errmessage{Use \string\frak\space only in math mode}}\fi\next}
\def\frak@#1{{\fam\eufam{{#1}}}}
\catcode`@=12

\newcommand{\be}{\begin{equation}}
\newcommand{\ee}{\end{equation}}

\newcommand{\bea}{\begin{eqnarray}}
\newcommand{\eea}{\end{eqnarray}}

\newcommand{\eq}{\! & = & \!}
\newcommand{\eqq}{\!\! & = & \!\!}

\newcommand{\hst}[1]{\rule{#1}{0mm}}
 
\newcommand{\what}[1]{\widehat{\rule{#1}{0pt}}}

\newcommand{\nn}{\nonumber}

\newcommand{\sect}[1]{\setcounter{equation}{0}\section{#1}}

\newcommand{\reff}[1]{(\ref{#1})}

\newcommand{\Ga}{\Gamma}
\newcommand{\De}{\Delta}
\newcommand{\Lam}{\Lambda}
\newcommand{\al}{\alpha}
\newcommand{\bet}{\beta}
\newcommand{\ga}{\gamma}
\newcommand{\de}{\delta}
\newcommand{\ep}{\epsilon}
\newcommand{\ve}{\varepsilon}
\newcommand{\io}{\iota}
\newcommand{\lam}{\lambda}
\newcommand{\vr}{\varrho}
\newcommand{\si}{\sigma}
\newcommand{\vp}{\varphi}
\newcommand{\om}{\omega}

\newcommand{\la}{\langle}
\newcommand{\ra}{\rangle}
\newcommand{\rar}{\rightarrow}
\newcommand{\lra}{\longrightarrow}
\newcommand{\ti}{\times}
\newcommand{\op}{\oplus}
\newcommand{\bop}{\bigoplus}
\newcommand{\ot}{\otimes}
\newcommand{\otb}{{\bar{\otimes}}}
\newcommand{\otbep}%
{\bar{\otimes}_{\rule{-0.05em}{0em}{}_{{\scriptstyle\varepsilon}}}}
\newcommand{\otbeppr}%
{\bar{\otimes}_{\rule{-0.1em}{0em}{}_{{\scriptstyle\varepsilon}'}}}
\newcommand{\pa}{\partial}

\newcommand{\KK}{\Bbb K}
\newcommand{\ZZ}{\Bbb Z}
\newcommand{\NN}{\Bbb N}
\newcommand{\NS}{{\Bbb N}_{\ast}}

\newcommand{\Ze}{\{0\}}
\newcommand{\Lg}{{\mbox{Lgr}}}
\newcommand{\Lgn}{{\mbox{Lgr}_n}}
\newcommand{\Lgnpo}{{\mbox{Lgr}_{n+1}}}
\newcommand{\cd}{\cdot}
\newcommand{\od}{\odot}
\newcommand{\Sn}{{{\cal S}_n}}
\newcommand{\ad}{{\mbox{ad}}}
\newcommand{\Hh}{\hat{H}}
\newcommand{\Hhastgr}{\hat{H}^{\ast\mbox{{\scriptsize gr}}}}
\newcommand{\Lh}{\hat{L}}
\newcommand{\pih}{\hat{\pi}}
\newcommand{\Et}{\tilde{E}}
\newcommand{\Ht}{\tilde{H}}
\newcommand{\Lzb}{L_{\bar{0}}}
\newcommand{\Lob}{L_{\bar{1}}}

\newcommand{\Qpm}{Q_{\pm}}
\newcommand{\Upm}{U_{\pm}}
\newcommand{\Vpm}{V_{\pm}}
\newcommand{\Wpm}{W_{\pm}}
\newcommand{\Xpm}{X_{\pm}}
\newcommand{\epm}{e_{\pm}}
\newcommand{\vpm}{v_{\pm}}
\newcommand{\wpm}{w_{\pm}}
\newcommand{\hal}{\frac{1}{2}}
\newcommand{\ohal}{{\textstyle\frac{1}{2}}}
\newcommand{\thhal}{{\textstyle\frac{3}{2}}}
\newcommand{\nhal}{{\textstyle\frac{n}{2}}}
\newcommand{\zb}{\bar{0}}
\newcommand{\ob}{\bar{1}}
\newcommand{\tpb}{\overline{2p}}
\newcommand{\tqb}{\overline{2q}}
\newcommand{\Vfb}{\bar{V}_4}

\newcommand{\comp}{\!\stackrel{\textstyle{\hst{0ex} \atop \circ}}{\hst{0ex}}\!}
\newcommand{\wh}{\widehat}

\newcommand{\wepf}{{\:\wedge\!\!\!{}_{{}_{{}_{\scriptstyle\varepsilon}}}\:}}
\newcommand{\wepe}{{\wedge\!\!\!{}_{{}_{{}_{\scriptstyle\varepsilon}}}}}
\newcommand{\Wepf}%
{{\:\bigwedge\!\!\!\!\!\!\:{}_{{}_{{}_{{}_{{}_{{}_{\scriptstyle%
\varepsilon}}}}}}\:}}
\newcommand{\Wepfz}%
{{\:\bigwedge^{0}\!\!\!\!\!\!\:{}_{{}_{{}_{{}_{{}_{{}_{\scriptstyle%
\varepsilon}}}}}}\:}}
\newcommand{\Wepfo}%
{{\:\bigwedge^{1}\!\!\!\!\!\!\:{}_{{}_{{}_{{}_{{}_{{}_{\scriptstyle%
\varepsilon}}}}}}\:}}
\newcommand{\Wepft}%
{{\:\bigwedge^{2}\!\!\!\!\!\!\:{}_{{}_{{}_{{}_{{}_{{}_{\scriptstyle%
\varepsilon}}}}}}\:}}
\newcommand{\Wepfd}%
{{\:\bigwedge^{3}\!\!\!\!\!\!\:{}_{{}_{{}_{{}_{{}_{{}_{\scriptstyle%
\varepsilon}}}}}}\:}}
\newcommand{\Wepfn}%
{{\:\bigwedge^{n}\!\!\!\!\!\!\:{}_{{}_{{}_{{}_{{}_{{}_{\scriptstyle%
\varepsilon}}}}}}\:}}
\newcommand{\Wepfnpe}%
{{\:\bigwedge^{n+1}\rule{-1em}{0ex}{}_{{}_{{}_{{}_{{}_{{}_{\scriptstyle%
\varepsilon}}}}}}\:}}
\newcommand{\Wept}%
{{\bigwedge\!\!\!\!{}_{{}_{{}_{{}_{{}_{\scriptstyle\varepsilon}}}}}\!\;}}
\newcommand{\Weptt}%
{{\bigwedge\limits^{2}\!\!\!\!{}_{{}_{{}_{{}_{{}_{\scriptstyle%
\varepsilon}}}}}\!\;}}
\newcommand{\Weptn}%
{{\bigwedge\limits^{n}\!\!\!\!{}_{{}_{{}_{{}_{{}_{\scriptstyle%
\varepsilon}}}}}\!\;}}

\newcommand{\ario}{{\,\stackrel{\iota}{\longrightarrow}\,}}
\newcommand{\ariopr}{{\,\stackrel{{\iota}'}{\longrightarrow}\,}}
\newcommand{\ariob}{{\,\stackrel{\bar{\iota}}{\longrightarrow}\,}}
\newcommand{\arioh}{{\,\stackrel{\hat{\iota}}{\longrightarrow}\,}}
\newcommand{\arpi}{{\,\stackrel{\pi}{\longrightarrow}\,}}
\newcommand{\arpipr}{{\,\stackrel{{\pi}'}{\longrightarrow}\,}}
\newcommand{\arpib}{{\,\stackrel{\bar{\pi}}{\longrightarrow}\,}}
\newcommand{\arpih}{{\,\stackrel{\hat{\pi}}{\longrightarrow}\,}}
\newcommand{\arvp}{{\,\stackrel{\varphi}{\longrightarrow}\,}}
\newcommand{\arid}{{\,\stackrel{id}{\longrightarrow}\,}}
\newcommand{\arnu}{{\,\stackrel{\nu}{\longrightarrow}\,}}

\begin{document}

\begin{titlepage}

\hspace*{\fill} \parbox{9em}{BONN--TH--97--01 \\
                             January 1997        }

\vspace{18mm}

\begin{center}
{\LARGE\bf Cohomology of Lie superalgebras \\[0.5ex] 
           and of their generalizations} \\
\vspace{11mm}
{\large M. Scheunert} \\
\vspace{1mm}
Physikalisches Institut der Universit\"{a}t Bonn \\
Nu{\ss}allee 12, D--53115 Bonn, Germany \\
\vspace{5mm}
{\large R.B. Zhang} \\
\vspace{1mm}
Department of Pure Mathematics \\
University of Adelaide \\ 
Adelaide, Australia \\
\end{center}

\vspace{15mm}
\begin{abstract}
\noindent
The cohomology groups of Lie superalgebras and, more generally, of $\ve$
Lie algebras, are introduced and investigated. The main emphasis is on the
case where the module of coefficients is non--trivial. Two general
propositions are proved, which help to calculate the cohomology groups.
Several examples are included to show the peculiarities of the super case. For
$L = sl(1|2)$, the cohomology groups $H^1(L,V)$ and $H^2(L,V)$, with $V$ a
finite--dimensional simple graded $L$--module, are determined, and the result
is used to show that $H^2(L,U(L))$ (with $U(L)$ the enveloping algebra of $L$)
is trivial. This implies that the superalgebra $U(L)$ does not admit of any
non--trivial formal deformations (in the sense of Gerstenhaber). Garland's
theory of universal central extensions of Lie algebras is generalized to the
case of $\ve$ Lie algebras.
\end{abstract}

\vspace{\fill}

q-alg/9701037

\end{titlepage}

\setcounter{page}{2}
\renewcommand{\baselinestretch}{1}     
\small\normalsize                      
\noindent
{\Large \bf
I. Introduction \\[-0.5ex]} \setcounter{section}{1}

\noindent
Cohomology is an important tool in mathematics. Its range of applications
contains algebra and topology as well as the theory of smooth manifolds
or of holomorphic functions. The cohomology theory of Lie algebras (whose
generalization is the main objective of the present paper) has its origins
in the work by E. Cartan, but the foundation of the theory, as an
independent topic of research, is due to Chevalley, Eilenberg \cite{ChE},
to Koszul \cite{Ko}, and to Hochschild, Serre \cite{HSe}. A unifying
treatment of the cohomology theory of groups, associative algebras, and
Lie algebras has been given by Cartan, Eilenberg \cite{CaE}.

As is well--known, several classical results in Lie algebra theory have a
cohomological interpretation. Thus, if $L$ is a Lie algebra, the structure of
the extensions of $L$--modules (and hence the question of semi--simplicity
of $L$--modules) is described by the 1--cohomology of $L$\,, and the
structure of the Lie algebra extensions (and hence the Levi, Malcev theorem)
is related to the 2--cohomology \cite{ChE,SSL}.

In view of these remarks, it is hardly surprising that shortly after the
birth of supersymmetry, a paper appeared \cite{Le} in which the most basic
constructions and results of the classical theory are generalized to the
case of Lie superalgebras \cite{Ka,S}. Later on, Leites and Fuks calculated
the cohomology groups of the classical Lie superalgebras with trivial
coefficients \cite{FL} (see also Ref.~\cite{Fu}), and Leites and V. Serganova
used cohomological methods to determine systems of generators and relations
for these algebras \cite{LS}. In addition, we should mention a letter by
Retakh and Feigin \cite{RF} and the papers by Tripathy and coworkers (see
Ref.~\cite{TP} and the papers cited therein). Apart from this, very little
seems to be known about the cohomology of Lie superalgebras, in particular,
about the cohomology with non--trivial coefficients.

After the present work was put into the internet and submitted for publication,
C. Gruson sent us a copy of the galley proofs of an article \cite{Gr}, in which
she proves the finiteness of the homology of certain modules over a Lie
superalgebra. She also drew our attention to a paper by J. Tanaka \cite{Ta},
in which the author calculates the dimensions of the homology and cohomology
groups of $sl(2|1)$ with coefficients in an arbitrary finite--dimensional
simple module. This has a direct bearing on part of our work, and we shall
comment on it later on.

Our present interest in this topic has been stimulated by the theory of
quantum algebras \cite{Dr,Ji}. It is known that, for any semi--simple Lie
algebra $L$\,, the enveloping algebra $U(L)$ does not admit of any
non--trivial formal deformations in the sense of Gerstenhaber \cite{Ge}.
We would like to know for which of the simple Lie superalgebras the same
holds true. Proceeding as in Ref.~\cite{Kas}, this amounts to showing that
the second cohomology group $H^2(L,U(L))$ is trivial. Unfortunately, we
are far from proving a general result of this type. In fact, all we can say
at present is that this is true if $L$ is one of the algebras $osp(1|2n)$
(this is easy to see and, certainly, was already known to Kac \cite{Ka}), or
if $L$ is equal to $sl(1|2)$ (this will be proved in the present work).
Moreover, in the meantime we have proved that $H^2(L,U(L))$ is also trivial
for the algebras $sl(m|1)$ with $m \geq 3$\,. This will be published elsewhere.
Accordingly, our paper is going to be explorative in character.

Let us briefly describe the setup of the present article. In Sec.~II we
collect the basic definitions and constructions, which will be needed in
the subsequent sections. Using the graded multilinear algebra as described
in Ref.~\cite{S.ma} this is just a simple transcription of the classical
situation to the graded case. In doing so, we shall follow the approach of
Ref.~\cite{HSe} (see also Ref.~\cite{ChE}). In addition, we comment on some
peculiarities of the graded case, and we prove two propositions, which will
be the main tools in our calculations. The first one shows that it is
sufficient to investigate cochains and coboundaries with certain invariance
properties, the second proposition gives a general criterion for the
cohomology with values in a given module to be trivial.

Sec.~III contains several examples. In the first example we consider the
cohomology of semi--simple Lie algebras. This example is included mainly in
order to show how far we can get, under the most favourable circumstances,
by only using the results of Sec.~II. A short look at this example reveals
that the same arguments can be applied to the Lie superalgebras $osp(1|2n)$.
This will be the content of the second example. The third example contains
some special information on the cohomology with trivial coefficients (again
well--known in the case of semi--simple Lie algebras). In the fourth example
we use the second proposition of Sec.~II to determine, for the general linear
Lie superalgebras $gl(m|n)$, which of the finite--dimensional simple graded
$gl(m|n)$--modules might have a non--trivial cohomology. This discussion
also covers the case of the special linear Lie superalgebras $sl(m|n)$,
provided that $m \neq n$\,.

In Sec.~IV we generalize Garland's theory of universal central coverings
(i.e., extensions) of Lie algebras \cite{Ga} to the general graded setting.
As a simple example, we determine the universal covering of the simple Lie
superalgebras $sl(n|n)/\mbox{center}$, with $n \ge 2$.

In Sec.~V we consider the Lie superalgebra $L = sl(1|2)$ and determine the
finite--dimensional simple graded $L$--modules $V$ for which $H^1(L,V)$ or
$H^2(L,V)$ is non--trivial. The result is then used to show that
$H^2(L,U(L))$ is trivial. Of course, the modules $V$ in question can be read
off from Ref.~\cite{Ta}. Nevertheless, we think it is worthwhile to present
our own calculations. First of all, this makes our paper self--contained in
this respect. Secondly, since we only determine $H^1(L,V)$ and $H^2(L,V)$,
our calculations are much simpler and the results much more explicit than
those of Ref.~\cite{Ta}, where $\dim H^n(L,V)$ is determined for arbitrary
$n$\,. Finally, this part of our paper is meant to explain how one might
investigate $H^2(L,U(L))$ for more general Lie superalgebras $L$\,.

Sec.~VI contains a brief discussion of our results. The paper is closed by
two appendices. In App.~A we describe the analogue of the exterior (i.e.,
Grassmann) algebra in the general graded setting. App.~B contains some
special information related to Sec.~V. This will be used to show that, for
$L = sl(1|2)$, the cohomology group $H^1(L,U(L))$ is non--trivial.

We close this introduction by explaining some of our conventions. The base
field will always be a field $\KK$ of characteristic zero. We use the
standard notation $\ZZ$ for the ring of integers and $\NN$ for the set of
natural numbers (according to our convention, $\NN$ contains $0$). The set of
strictly positive integers will be denoted by $\NS$\,. As explained above, in
Sec.~II and Sec.~IV as well as in App.~A we are going to consider the so--called
$\ve$ Lie algebras \cite{S.gl} (also called colour Lie algebras \cite{RW}).
In these sections, $\Ga$ is an arbitrary Abelian group and $\ve$ is a
commutation factor on $\Ga$ with values in $\KK$\,. The reader who is only
interested in Lie superalgebras may simply set $\Ga = \ZZ_2$ or $\Ga = \ZZ$
and choose $\ve$ to be the standard commutation factor of supersymmetry,
defined by $\ve(\al,\bet) = (-1)^{\al\bet}$ for all $\al,\bet \in \Ga$.
The choice $\Ga = \ZZ_2$ leads to general Lie superalgebras, for $\Ga = \ZZ$\,,
we obtain the consistently $\ZZ$--graded Lie superalgebras in the sense of
Ref.~\cite{Ka}. The multiplication in an $\ve$ Lie algebra (and hence in a
Lie superalgebra) will be denoted by a pointed bracket $\la\;,\;\ra$.
\vspace{5ex}

%
\noindent
{\Large \bf
\begin{tabbing}
II.\ \= \kill
II.  \> Cohomology of $\ve$ Lie algebras: \\[0.5ex]
     \> Definitions and some basic results
\end{tabbing}\vspace{0ex}} \setcounter{section}{2} \setcounter{equation}{0}

\noindent
Let us first introduce the notation which will be used throughout the present
section. Quite generally, we adopt the conventions of Ref.~\cite{S.ma}. In
the following, $\Ga$ denotes an Abelian group, $\ve : \Ga \ti \Ga \rar \KK$
a commutation factor on $\Ga$ with values in $\KK$\,, $L$ an $\ve$ Lie
algebra \cite{S.gl}, and $V$ a graded $L$--module. All gradations are
understood to be $\Ga$--gradations. If an element of a graded vector space is
homogeneous of a certain degree, then this degree will be denoted by the
lower case ``Greek analogue'' of the letter denoting the element itself. For
example, if $A, b', g, x$ are homogeneous elements, their degrees will be
denoted by $\al, \bet', \ga, \xi \in \Ga$, respectively. Only occasionally
(but in all cases where there could be any doubts) the degrees will be
specified explicitly.

Let $n \ge 1$ be a natural number and let
\be  \Lgn(L,V) = \Lgn(\underbrace{L,\ldots,L}_{n};V) \ee
be the graded vector space of all $n$--linear mappings of
$L^n = L \ti\ldots\ti L$ into $V$ which can be written as a sum of homogeneous
$n$--linear mappings of $L^n$ into $V$. Recall that $\Lgn(L,V)$ is equal to
the space of all $n$--linear mappings of $L^n$ into $V$ if (for example)
$\Ga$ is finite or if $L$ is finite--dimensional.

According to Ref.~\cite{S.ma} $\Lgn(L,V)$ is a graded $L$--module: If
$A \in L$ and $g \in \Lgn(L,V)$ are homogeneous, the action of $A$ on $g$ is
defined by
\newpage  
\bea (A_{L}(g))(A_1,\ldots,A_n) \eq A_{V}(g(A_1,\ldots,A_n)) \\
   & & \!\!\mbox{}-\sum_{r=1}^{n}\ve(\al,\ga + \al_1 + \ldots + \al_{r-1})\,
       g(A_1,\ldots,\la A,A_r\ra,\ldots,A_n) \,, \nn
\eea
for all homogeneous elements $A_1,\ldots,A_n \in L$\,. (Recall that the
$\ve$--commutator of two elements $A,B \in L$ is denoted
by $\la A,B \ra$. Quite generally, if $V$ is a graded $L$--module, the
representative of an element $A \in L$ in $V$ will be denoted by $A_V$\,.
If $x$ is an element of $V$, we shall frequently simplify the notation and
write $A \cd x$ in place of $A_V(x)$.)

On the other hand, again according to Ref.~\cite{S.ma}, there is a natural
representation of the symmetric group $\Sn$ in $\Lgn(L,V)$. To describe it,
we define a map
\be  \ve_n : \Sn \ti \Ga^n \lra \KK   \ee
by
\be  \ve_n(\pi^{-1};\ga_1,\ldots,\ga_n)\; =
               \prod_{i < j \atop \pi(i) > \pi(j)}\!\!\ve(\ga_i,\ga_j) \ee
for all $\pi \in \Sn$ and $\ga_1,\ldots,\ga_n \in \Ga$. The basic property
of this map is that
\be  \ve_n(\pi\tau;\ga_1,\ldots,\ga_n) =
\ve_n(\pi;\ga_1,\ldots,\ga_n)\,\ve_n(\tau;\ga_{\pi(1)},\ldots,\ga_{\pi(n)}) \ee
for all $\pi,\tau \in \Sn$ and all $\ga_1,\ldots,\ga_n \in \Ga$. Then the
representation $\pi \rar \check{S}_{\pi}$ of $\Sn$ in $\Lgn(L,V)$ is given by
\be (\check{S}_{\pi}(g))(A_1,\ldots,A_n)
         = \ve_n(\pi;\al_1,\ldots,\al_n)\,g(A_{\pi(1)},\ldots,A_{\pi(n)})
                                                       \label{checkS} \ee
for all $\pi \in \Sn$\,, $g \in \Lgn(L,V)$, and all homogeneous elements
$A_1,\ldots,A_n \in L$\,. It is known that $\check{S}_{\pi}$ is an automorphism of the graded $L$--module
$\Lgn(L,V)$, for all $\pi \in \Sn$\,.

Now let $C^n(L,V)$ be the set of all $\ve$--skew--symmetric elements of
$\Lgn(L,V)$, i.e., the set of all elements $g \in \Lgn(L,V)$ such that
\be  \check{S}_{\pi}\,g = \mbox{sgn}(\pi)\,g   \ee
for all $\pi \in \Sn$\,, where $\mbox{sgn}(\pi)$ denotes the signum of the
permutation $\pi$\,. We extend this definition to the case of integers
$n \le 0$ and set
\bea       C^n(L,V) \eq \Ze \;\;\;\mbox{if}\;\;\; n \le -1 \\
           C^0(L,V) \eq V \,. \label{cochz}
\eea
Occasionally, it is useful to identify $C^0(L,V)$ with $\Lg(\KK\,,V)$. Then
$C^n(L,V)$ can be identified with $\Lg(\Weptn L,V)$, where $\Weptn L$ is the
subspace of degree $n$ of the $\ve$--exterior algebra $\Wept L$ of $L$\,. A few
comments on the $\ve$--exterior algebra of a $\Ga$--graded vector space are
contained in App.~A. Obviously, $C^n(L,V)$ is a graded $L$--module, for all~%
$n$\,. The elements of $C^n(L,V)$ are called the $n${\em --cochains} of $L$
with values in $V$.

In the subsequent discussion we follow Ref.~\cite{HSe}. First we define, for
all integers $n$\,, a map
\be  C^n(L,V) \ti L \lra C^{n-1}(L,V) \,,  \ee
denoted by
\be  (g,A) \lra g_A \,,  \ee
by
\be  g_A = 0 \;\;\;\mbox{if}\;\;\; n \le 0 \,,  \ee
\be  g_A = g(A) \;\;\;\mbox{if}\;\;\; n = 1 \,,  \ee
\be  g_A(A_2,\ldots,A_n) = g(A,A_2,\ldots,A_n) \;\;\;\mbox{if}\;\;\;
                                                             n \ge 2 \,, \ee
where $g \in C^n(L,V)$ and $A,A_2,\ldots,A_n \in L$\,. It is easy to see that
this map is bilinear, homogeneous of degree zero, and $L$--invariant, i.e.,
\be  B \cd g_A = (B \cd g)_A + \ve(\bet,\ga)\,g_{\la B,A\ra}  \ee
for all homogeneous elements $A,B \in L$ and $g \in C^n(L,V)$.

Next we define, again for all integers $n$\,, the linear {\em coboundary
operator}
\be  \de^n : C^n(L,V) \lra C^{n+1}(L,V) \label{del}  \ee
as follows: We set
\be  \de^n = 0 \;\;\;\mbox{if}\;\;\; n \le -1  \ee
and define $\de^n$ for $n \ge 0$ inductively by
\be  (\de^n(g))_A = \ve(\ga,\al)A \cd g - \de^{n-1}(g_A) \,,\label{delind} \ee
where $g \in C^n(L,V)$ and $A \in L$ are homogeneous. (Note that
Eq.~\reff{delind} is trivially satisfied if $n \le -1$\,.) As it stands,
Eq.~\reff{delind} defines a linear map of $C^n(L,V)$ into $\Lgnpo(L,V)$ which
is homogeneous of degree zero, but it is easy to see by induction on $n \ge 0$
that $\de^{n}g$ is $\ve$--skew--symmetric and hence an element of
$C^{n+1}(L,V)$. One can now show by induction on $n$ that $\de^n$ is a
homomorphism of graded $L$--modules, i.e., that
\be  A \cd(\de^n g) = \de^n(A \cd g) \label{delinv}  \ee
for all $n$ and all $A \in L$\,, $g \in C^n(L,V)$, and then (very easily)
that
\be  \de^{n+1}\comp\de^n = 0 \label{deldel} \ee
for all $n$\,. Finally, one can verify inductively that the following
explicit formula holds for $n \ge 1$\,:
\bea \lefteqn{(\de^n g)(A_0,A_1,\ldots,A_n)} \hspace{-0.2em} \label{delex} \\
 \eq \sum_{r=0}^{n}(-1)^r\ve(\ga+\al_0+\ldots+\al_{r-1},\al_r)
                          \,A_r \cd g(A_0,\ldots,\wh{A_r},\ldots,A_n) \nn \\
 & & \hspace{-1.4em}\mbox{}+\sum_{r < s}(-1)^s
 \ve(\al_{r+1}+\ldots+\al_{s-1},\al_s)\,
 g(A_0,\ldots,A_{r-1},\la A_r,A_s\ra,A_{r+1},\ldots,\wh{A_s},\ldots,A_n)\,,\nn
\eea
where $g \in C^n(L,V)$ and $A_0,A_1,\ldots,A_n \in L$ are homogeneous, and
where the sign $\what{0.5em}$ indicates that the element below it must be
omitted. Empty sums (like $\al_0+\ldots+\al_{r-1}$ for $r = 0$ and
$\al_{r+1}+\ldots+\al_{s-1}$ for $s = r+1$) are set equal to zero. In
particular, for $n = 1$ we have
\be  (\de^1 g)(A_0,A_1) = \ve(\ga,\al_0)\,A_0\cd g(A_1)
 -\ve(\ga+\al_0,\al_1)\,A_1\cd g(A_0) - g(\la A_0,A_1\ra) \,,\label{del1} \ee
and for $n = 2$ we obtain
\bea (\de^2 g)(A_0,A_1,A_2)
 \eqq \ve(\ga,\al_0)\,A_0\cd g(A_1,A_2)
               - \ve(\ga+\al_0,\al_1)\,A_1\cd g(A_0,A_2) \label{del2} \\
 & & \hspace{-0.6em}\mbox{}+\ve(\ga+\al_0+\al_1,\al_2)\,A_2\cd g(A_0,A_1) \nn\\
 & & \hspace{-0.6em}\mbox{}-g(\la A_0,A_1\ra,A_2)
    + \ve(\al_1,\al_2)\,g(\la A_0,A_2\ra,A_1) + g(A_0,\la A_1,A_2\ra) \,. \nn
\eea
Occasionally, it is useful to rewrite Eq.~\reff{delex} in the form 
\bea \lefteqn{2(\de^n g)(A_0,A_1,\ldots,A_n)} \hspace{2em} \label{delalt} \\
 \eq \sum_{r=0}^{n}(-1)^r\ve(\ga+\al_0+\ldots+\al_{r-1},\al_r)
                         (A_r \cd g)(A_0,\ldots,\wh{A_r},\ldots,A_n) \nn \\
 & & \hspace{-1.4em}\mbox{}
         + \sum_{r=0}^{n}(-1)^r\ve(\ga+\al_0+\ldots+\al_{r-1},\al_r)
                         \,A_r \cd g(A_0,\ldots,\wh{A_r},\ldots,A_n) \,. \nn
\eea

\noindent
{\em Remark 2.1.} The relations \reff{delind}, \reff{delinv}, \reff{deldel}
are the basic properties of the cohomology operators $\de^n$. Of course, they
can be checked directly (regarding Eq.~\reff{delex} as the definition of
$\de^n$\,), but it is much easier to proceed as in Ref.~\cite{HSe}, i.e., to
follow the steps described above.

\noindent
{\em Remark 2.2.} We could use the $\ve$--skew--symmetry of $g$ to shift the
argument $\la A_r,A_s\ra$ in the second sum on the right hand side of
Eq.~\reff{delex} to the first place. Then this equation would take a form
which is similar to the one familiar from the non--graded theory. However,
this would lead to additional $\ve$ factors which we wanted to avoid.

Next let $Z^n(L,V)$ denote the kernel of $\de^n$ and let $B^n(L,V)$ denote
the image of $\de^{n-1}$. Eq.~\reff{delinv} implies that $Z^n(L,V)$ and
$B^n(L,V)$ are graded submodules of $C^n(L,V)$, and according to
Eq.~\reff{deldel} we have
\be  B^n(L,V) \subset Z^n(L,V) \,. \label{BinZ} \ee
The elements of $Z^n(L,V)$ are called $n${\em --cocycles}, the elements of
$B^n(L,V)$ are the $n${\em --coboundaries.} Because of \reff{BinZ} we can
construct the so--called {\em cohomology groups}
\be  H^n(L,V) = Z^n(L,V)/B^n(L,V) \,. \ee
Two elements of $Z^n(L,V)$ are said to be {\em cohomologous} if their
residue classes modulo $B^n(L,V)$ coincide, i.e., if their difference lies
in $B^n(L,V)$.

Of course, the cohomology groups are graded $L$--modules, too. Actually,
their $L$--module structure is trivial. To prove this we have to show that
for any cocycle $g \in Z^n(L,V)$ and any $A \in L$\,, the cocycle $A \cd g$
is a coboundary, and this follows directly from Eq.~\reff{delind}.

The following consequence of the foregoing result can be used to facilitate
the calculation of the cohomology groups.

\noindent
{\bf Proposition 2.1}\\
{\it We use the notation introduced at the beginning of the present section.
Let $L'$ be a graded subalgebra of $L$\,.\\
a) If the $L'$--module $Z^{n}(L,V)$ is semi--simple (which is true if the
$L'$--module $C^n(L,V)$ is semi--simple), then for any cocycle $g \in Z^n(L,V)$
there exists an $L'$--invariant cocycle $g' \in Z^n(L,V)$ which is
cohomologous to $g$\,. \\
b) If the $L'$--module $C^{n-1}(L,V)$ is semi--simple, then any $L'$--invariant
coboundary $b \in B^n(L,V)$ is equal to $\de^{n-1}(g)$ with an $L'$--invariant
cochain $g \in C^{n-1}(L,V)$. \\
(Semi--simplicity may be understood in the graded sense.)}

\noindent
Proof \\
a) If the graded $L'$--module $Z^n(L,V)$ is semi--simple, there exists a
graded $L'$--submodule $X^n(L,V)$ of $Z^n(L,V)$ which is complementary to
$B^n(L,V)$. By definition, the canonical map
\be  Z^n(L,V) \lra H^n(L,V) \ee
induces an isomorphism
\be  X^n(L,V) \lra H^n(L,V) \ee
of graded $L'$--modules. Since the $L'$--module structure of $H^n(L,V)$ is
trivial, the same is true for the $L'$--module structure of $X^n(L,V)$, i.e.,
all elements of $X^n(L,V)$ are $L'$--invariant. Thus all we have to do is to
choose an element $g' \in X^n(L,V)$ such that the canonical images of $g$
and $g'$ in $H^n(L,V)$ coincide.

\noindent
b) If the graded $L'$--module $C^{n-1}(L,V)$ is semi--simple, there exists a
graded $L'$--submodule $Y^{n-1}(L,V)$ of $C^{n-1}(L,V)$ which is
complementary to the kernel of $\de^{n-1}$ (i.e., to $Z^{n-1}(L,V)$). Then
$\de^{n-1}$ induces an isomorphism
\be  Y^{n-1}(L,V) \lra B^n(L,V) \ee
of graded $L'$--modules. This implies b) and completes the proof of the
proposition.

Next we consider some constructions which can be carried out with the
cohomology groups. Let $L'$ be a second $\ve$ Lie algebra and let
\be  \om : L' \lra L  \ee
be a homomorphism of $\ve$ Lie algebras. Then the action of $L'$ on $V$,
defined by
\be  A' \cd x = \om(A') \cd x  \ee
for all $A' \in L'$ and $x \in V$, makes $V$ into a graded $L'$--module,
which we denote by $V^{\om}$. If $g \in C^n(L,V)$ with some $n \ge 1$\,,
then the map $L'^{\,n} \rar V$, defined by
\be  (A'_1,\ldots,A'_n) \lra g(\om(A'_1),\ldots,\om(A'_n))  \ee
for all $A'_1,\ldots,A'_n \in L'$, is an element of $C^n(L',V^{\om})$. This
assignment defines a linear map
\be  C^n(L,V) \lra C^n(L',V^{\om}) \,. \label{omchain} \ee
We extend this definition to the cases $n \le 0$ by defining this map to be
equal to $id_V$ if $n = 0$\,, and to zero if $n \le -1$\,. Obviously, the
map \reff{omchain} is homogeneous of degree zero, and it is easy to see that
it is compatible with the coboundary operators, in the obvious sense. In
particular, it maps cocycles into cocycles and coboundaries into coboundaries,
thus inducing a corresponding map for the cohomology groups:
\be H^n(L,V) \lra H^n(L',V^{\om}) \,. \label{omcohom} \ee

Similarly, let
\be  f : V \lra W  \ee
be a linear map of $V$ into a second graded $L$--module $W$ and suppose that
$f$ is homogeneous of degree $\vp$ and $L$--invariant. Recall \cite{S.ma}
that $L$--invariance means that
\be  f(A \cd x) = \ve(\vp,\al)\,A \cd f(x)  \ee
for all $x \in V$, and for all $A \in L$ which are homogeneous of degree $\al$
(stated differently, $f$ is $L$--linear in the graded sense). For any
integer $n$\,, we define the map
\be  f^n_c : C^n(L,V) \lra C^n(L,W)  \ee
by
\be  f^n_c = 0 \;\;\;\mbox{if}\;\;\; n \le -1  \ee
\be  f^0_c = f \;\;\;\mbox{if}\;\;\; n = 0  \ee
\be  f^n_c(g) = f \comp g \;\;\;\mbox{if}\;\;\; n \ge 1 \,. \ee
Obviously, $f^n_c$ is homogeneous of degree $\vp$\,, and it is easy to see that
\be  \de^n_W \comp f^n_c = f^{n+1}_c \comp \de^n_V \,, \label{delmap}  \ee
where, temporarily, we have marked the coboundary operators for cochains with
values in $V$ and $W$ by the subscripts $V$ and $W$, respectively.
Eq.~\reff{delmap} implies that $f^n_c$ maps cocycles into cocycles and
coboundaries into coboundaries, and hence induces a linear map
\be  f^n : H^n(L,V) \lra H^n(L,W) \,. \ee
Of course, $f^n$ is also homogeneous of degree $\vp$\,. If $U$ is a third
graded $L$--module and if $f': U \rar V$ is a linear map which is homogeneous
and $L$--invariant, then, obviously,
\be  (f \comp f')^n_c = f^n_c \comp f'^{n}_c \;\;\;\mbox{and}\;\;\;
                                   (f \comp f')^n = f^n \comp f'^{n} \,. \ee

As a first application of this construction, choose $\si \in \Ga$ and
consider the graded $L$--module $V^{\si}$ which is obtained from $V$ by a
shift of the gradation by $\si$ (see Ref.~\cite{S.ma}). Recall that,
considered as a vector space, $V$ and $V^{\si}$ coincide, and that the actions
of $L$ on $V$ and $V^{\si}$ are also the same, but that the gradations are
related by
\be  V^{\si}_{\ga} = V_{\ga + \si} \label{grshift} \ee
for all $\ga \in \Ga$. Define the linear map
\be  f : V^{\si} \lra V  \ee
by
\be  f(x) = \ve(\si,\xi)\,x  \ee
for all $x \in V^{\si}_{\xi}$\,, $\xi \in \Ga$. Then $f$ is bijective and
homogeneous of degree $\si$\,, and it is easy to see that $f$ is
$L$--invariant. Consequently,
\be  f^n : H^n(L,V^{\si}) \lra H^n(L,V) \label{hshift} \ee
is a bijective linear map which is homogeneous of degree $\si$\,. Stated
differently, $f^n$ is an isomorphism of the graded vector space
$H^n(L,V^{\si})$ onto the graded vector space $H^n(L,V)^{\si}$.

For our second application we suppose that $V$ is the direct sum of a family
$(V_i)_{i \in I}$ of graded submodules,
\be  V = \bigoplus_{i \in I}V_i \,.  \ee
Considering the canonical injections $V_i \rar V$ and projections
$V \rar V_i$\,, it is easy to construct a natural linear map
\be  \bigoplus_{i \in I}H^n(L,V_i) \lra H^n(L,V)  \label{dirsum}  \ee
which is injective and homogeneous of degree zero. If $I$ is finite or if $L$
is finite--dimensional, this map is even bijective, i.e., an isomorphism of
graded vector spaces. The details are obvious and may be left to the reader.

The third application is less trivial. Consider three graded $L$--modules
$U$, $V$, $W$ and a short exact sequence
\be \Ze \lra U \stackrel{j}{\lra} V \stackrel{p}{\lra} W \lra \Ze \,, \ee
where the linear mappings $j$ and $p$ are homogeneous (of degrees $\io$ and
$\pi$\,, respectively) and $L$--invariant. Then there exists a long exact
sequence
\be \ldots\, \lra H^n(L,U) \stackrel{j^n}{\lra} H^n(L,V)
           \stackrel{p^n}{\lra} H^n(L,W) \stackrel{\de^n}{\lra}
                   H^{n+1}(L,U) \lra \,\ldots\;, \label{cohomseq} \ee
where the maps $j^n$ and $p^n$ have been defined above. The definition of the
so--called {\em connecting homomorphisms}\, $\de^n$ (not to be confounded with
the coboundary operators $\de^n$ of Eq.~\reff{del}) is standard but a little
involved (see Refs.~\cite{McL,BT}). Since, in the applications, this
definition generally is of little use, we don't give it here but only mention
that $\de^n$ is homogeneous of degree $-\io -\pi$\,.

Next we want to construct certain product maps between the spaces of cochains.
Let $V$ and $W$ be two graded $L$--modules and let $m$ and $n$ be two integers.
We are going to define a bilinear map
\be  C^m(L,V) \ti C^n(L,W) \lra C^{m+n}(L,V \ot W) \,, \label{prod}  \ee
denoted by
\be  (g,h) \lra g \od h \,,  \ee
as follows. If $m < 0$ or $n < 0$\,, this map is equal to zero. Now assume
that $m,n \ge 0$\,, and let $g \in C^m(L,V)$ and $h \in C^n(L,W)$ be
homogeneous of degrees $\ga$ and $\eta$\,, respectively. If $m = n = 0$\,,
we set
\be  g \od h = g \ot h \,,  \ee
if $m = 0$\,, $n > 0$\,, we set
\be  (g \od h)(A_1,\ldots,A_n) = g \ot h(A_1,\ldots,A_n)  \ee
for all $A_1,\ldots,A_n \in L$\,, if $m > 0$\,, $n = 0$\,, we set
\be  (g \od h)(A_1,\ldots,A_m)
                = \ve(\eta,\al_1+\ldots+\al_m)\,g(A_1,\ldots,A_m) \ot h  \ee
for all homogeneous elements $A_1,\ldots,A_m \in L$\,. Finally, suppose that
$m,n \ge 1$\,. In that case, we define the $(m+n)$--linear map $g \otb h$
(the graded tensor product of $g$ and~$h$) by 
\bea \lefteqn{(g \otb h)(A_1,\ldots,A_m,A_{m+1},\ldots,A_{m+n})}
                                                            \hspace{2em} \\
\eq \ve(\eta,\al_1+\ldots+\al_m)\,g(A_1,\ldots,A_m)
                                       \ot h(A_{m+1},\ldots,A_{m+n}) \nn\eea
for all homogeneous elements $A_1,\ldots,A_{m+n} \in L$ and set
\be  g \od h
 = \frac{1}{m!\,n!}\sum_{\pi \in {\cal S}_{m+n}}\!\mbox{sgn}(\pi)
                             \check{S}_{\pi}(g \otb h) \,, \label{prodall}  \ee
where $\check{S}_{\pi}$ has been defined in Eq.~\reff{checkS}.

\noindent
{\it Remark 2.3.} If $m,n \ge 1$\,, let ${\cal S}(m,n)$ be the set of all
permutations $\pi \in {\cal S}_{m+n}$ such that $\pi$ is increasing on
$\{1,2,\ldots,m\}$ and on $\{m+1,m+2,\ldots,m+n\}$ (the appropriate set of
shuffle permutations). Then Eq.~\reff{prodall} can also be written in the form
\be  g \od h
 = \sum_{\pi \in {\cal S}(m,n)}\!\mbox{sgn}(\pi)
                             \check{S}_{\pi}(g \otb h) \,. \label{prodsim} \ee

Obviously, in all cases $g \od h$ is an element of $C^{m+n}(L,V \ot W)$.
Moreover, it follows from the general theory of graded $L$--modules \cite{S.ma}
that the bilinear map \reff{prod} is homogeneous of degree zero and
$L$--invariant, i.e., that
\be  A \cd (g \od h) = (A \cd g) \od h + \ve(\al,\ga)\,g \od (A \cd h)
                                                       \label{prodinv}  \ee
for all homogeneous elements $A \in L$\,, $g \in C^m(L,V)$, and
$h \in C^n(L,W)$.

The product is also compatible with the other operations introduced above.
In fact, using \reff{prodsim} it is not difficult to show that for all
integers $m,n$ and all elements $A \in L_{\al}$\,, $g \in C^m(L,V)$, and
$h \in C^n(L,W)_{\eta}$
 \be  (g \od h)_A = \ve(\eta,\al)\,g_A \od h + (-1)^m g \od h_A \,.  \ee
Using this relation and induction on $m+n$ it is easy to rederive
Eq.~\reff{prodinv} and then, again by induction on $m+n$\,, to obtain the
following generalized Leibniz rule
 \be  \de^{m+n}(g \od h) = (\de^m g) \od h + (-1)^m g \od(\de^n h) \,, \ee
which holds for all integers $m,n$ and all $g \in C^m(L,V)$ and
$h \in C^n(L,W)$.  An immediate consequence of this equation is that the
product maps \reff{prod} induce analogous ones for the cohomology groups:
 \be H^m(L,V) \ti H^n(L,W) \lra H^{m+n}(L,V \ot W) \,.  \ee
As a final result about the product $\od$ we notice that it is associative
in the obvious sense: If $U$ is a third graded $L$--module, we have for all
integers $\ell,m,n$
 \be  f \od (g \od h) = (f \od g) \od h \,, \label{prodass} \ee
where $f \in C^{\ell}(L,U)$, $g \in C^m(L,V)$, and $h \in C^n(L,W)$ (of
course, it is understood that $U \ot (V \ot W)$ and $(U \ot V) \ot W$ are
canonically identified). The proof can be carried out by induction on
$\ell + m + n$\,.

In the foregoing discussion, the product of two cocycles $g \in C^m(L,V)$ and
$h \in C^n(L,W)$ takes its values in $V \ot W$. This is the generic case, to
which the other possibilities can be reduced. If $U$ is a third graded
$L$--module and if
 \be  b : V \ti W \lra U  \ee
is a homogeneous $L$--invariant bilinear map, we can use the corresponding
linear map
 \be  \tilde{b} : V \ot W \lra U  \ee
to construct the maps
 \be  \tilde{b}^n_c : C^n(L,V \ot W) \lra C^n(L,U)  \ee
which, when composed with the product maps \reff{prod}, yield the maps
 \be  C^m(L,V) \ti C^n(L,W) \lra C^{m+n}(L,U) \,.  \ee
Of course, these have properties analogous to those of $\od$, in particular,
they induce bilinear maps
 \be H^m(L,V) \ti H^n(L,W) \lra H^{m+n}(L,U)  \ee
which are homogeneous of the same degree as $b$\,.

We close this section by showing that for certain coefficient modules $V$
the cohomology groups $H^n(L,V)$ must vanish. Let $U(L)$ be the enveloping
algebra of $L$ and let $U_{+}(L)$ denote the ideal of $U(L)$ which is
generated by $L$\,. We say that an element $X$ of $U(L)$ (and, by abuse of
language, also the representative $X_V$ of $X$ in an $L$--module $V$)
{\em does not contain a constant term,} in case $X \in U_{+}(L)$.

\noindent
{\bf Proposition 2.2}\\
{\it Let $L$ be a finite--dimensional $\ve$ Lie algebra, let $V$ be a graded
$L$--module, and let $C$ be a homogeneous element of the $\ve$ center of
$U(L)$ (i.e., a homogeneous Casimir element of $L$). Suppose that $C$ does not
contain a constant term, and that the operator $C_V$ (the representative of
$C$ in $V$) is invertible. Then we have
 \be  H^n(L,V) = \Ze \;\;\;\mbox{\it if}\;\;\; n \neq 0 \,.  \ee  }

\noindent
Proof \\
Let $(E_i)_{i \in I}$ be a basis of $L$\,, with $E_i$ homogeneous of degree
$\ep_i$\,. We define for all $A \in L$ and $j \in I$
 \be  \la A,E_j\ra = \sum_{i \in I}\ad_{ij}(A)\,E_i \:,  \ee
i.e., $A \rar (\ad_{ij}(A))$ is the adjoint representation, written in
matrix form.

Let $L^{\ast}$ be the coadjoint $L$--module, i.e., the graded dual of $L$
endowed with the representation contragredient to the adjoint one, and let
$(E'_i)_{i \in I}$ be the basis of $L^{\ast}$ dual to $(E_i)_{i \in I}$\,,
defined by
 \be E'_i(E_j) = \de_{ij} \;\;\;\mbox{for all}\;\;\; i,j \in I \,. \ee
Then $E'_i$ is homogeneous of degree $-\ep_i$\,.

Now consider an $L$--invariant $r$--linear form $\phi$ on $L^{\ast}$ which
is homogeneous of degree $\eta$\,. Then
 \be  C(\phi) = \sum_{i_1,\ldots,\,i_r \in I}
                  \phi(E'_{i_1},\ldots,E'_{i_r})\, E_{i_r}\ldots E_{i_1} \ee
is a Casimir element of $L$ (i.e., an element of the $\ve$ center of
$U(L)$) which is homogeneous of degree $\eta$ (see Ref.~\cite{S.ce}). We
define the following elements $C_i(\phi)$\,; $i \in I$, of $U(L)$\,:
 \be  C_i(\phi) = \sum_{i_2,\ldots,\,i_r \in I}
           \phi(E'_i,E'_{i_2},\ldots,E'_{i_r})\, E_{i_r}\ldots E_{i_2} \,. \ee
Obviously, we have
 \be  C(\phi) = \sum_{i \in I} C_i(\phi)\,E_i \,,  \ee
moreover, $C_i(\phi)$ is homogeneous of degree $\eta - \ep_i$\,, and the
$L$--invariance of $\phi$ implies that
 \be  \la A,C_i(\phi)\ra = - \sum_{j \in I}\ve(\al,\eta - \ep_j)
                                               \,\ad_{ij}(A)\,C_j(\phi)  \ee
for all homogeneous elements $A \in L$ and all $i \in I$.

Now let $C$ be the Casimir element described in the proposition, and suppose
that $C$ is homogeneous of degree $\eta$\,. Then it is known \cite{S.ce}
that $C$ can be written as a sum of Casimir elements of the form $C(\phi)$
considered above. Consequently, there exist elements $C_i \in U(L)$\,;
$i \in I$, which are homogeneous of degree $\eta - \ep_i$ and satisfy
 \be  C = \sum_{i \in I} C_i \, E_i  \ee
as well as
 \be  \la A,C_i\ra = - \sum_{j \in I}\ve(\al,\eta - \ep_j)
                                                    \,\ad_{ij}(A)\,C_j  \ee
for all homogeneous elements $A \in L$ and all $i \in I$.

After these preparations the proof of the proposition is easy. For any
integer $n \ge 1$ we define a linear map
 \be  d_n : C^n(L,V) \lra C^{n-1}(L,V)  \ee
by
 \be  (d_n(g))(A_2,\ldots,A_n) = \sum_{i \in I}\ve(\ep_i,\ga)
                                       \,C_i \cd g(E_i,A_2,\ldots,A_n)  \ee
for all homogeneous elements $g \in C^n(L,V)$ and all $A_2,\ldots,A_n \in L$
(with the obvious interpretation if $n = 1$). Then one can check that
 \be  d_{n+1} \comp \de^n + \de^{n-1} \comp d_n = (C_V)^n_c \,,
                                                           \label{hom}  \ee
where $C_V$ is the representative of $C$ in $V$ (see the analogous but
simpler calculation in the proof of Prop.~VII.5.6 in Ref.~\cite{HSt}). Of
course, $C_V$ is homogeneous of degree $\eta$ and $L$--invariant. Thus if
$g \in C^n(L,V)$ is a cocycle, the Eqs.~\reff{hom} and \reff{delmap} imply that
 \be  g = \de^{n-1}((C_V^{-1})_c^{n-1}d_n(g)) \,,  \ee
hence $g$ is a coboundary and the proposition is proved.

The case $n = 0$ (not covered by the proposition) can easily be dealt with:
Since $B^0(L,V) = \Ze$, we have 
 \be  H^0(L,V) = Z^0(L,V) = \{x \in V \,|\, A \cd x = 0
                     \,\;\mbox{for all}\;\, A \in L \} \,, \label{hzero} \ee
i.e., $H^0(L,V)$ is the space of all $L$--invariant elements of $V$.

\noindent
{\it Remark 2.4.} The definitions and constructions of the present section
(but not the two propositions) can be generalized by choosing, for the basic
domain of scalars, instead of the field $\KK$ an arbitrary associative
$\ve$--commutative graded algebra $S$ over $\KK$\,. This follows immediately
from the multilinear algebra over $S$ as sketched in the appendix of
Ref.~\cite{S.p}. A generalization of this type could be useful if one
considers the deformations of a superalgebra over a Grassmann algebra, and
the deformation parameter is allowed to be odd. (The algebra $L$ may also act
on $S$ by $\ve$--derivations; for the super case, see Ref.~\cite{Le}).
\vspace{5ex}

%
\noindent
{\Large \bf
III. Examples \\[-0.5ex]} \setcounter{section}{3} \setcounter{equation}{0}

\noindent
{\bf 1. Semi--simple Lie algebras}\\[0.8ex]
In our first example we consider a (finite--dimensional) semi--simple Lie
algebra $L$ and a finite--dimensional $L$--module $V$, and we want to
describe the cohomology groups $H^n(L,V)$. Needless to say, this case is
well--known, we consider it here in order to show how our propositions can
be applied under the most favourable circumstances.

The main input is that all finite--dimensional $L$--modules are semi--simple
(i.e., the corresponding representations are completely reducible). Because
of the isomorphism \reff{dirsum} we may, therefore, assume that the
$L$--module $V$ is simple (of course, in practice it may be difficult to
decompose an $L$--module into simple submodules). If $V$ is non--trivial, it
is well--known that there exists a quadratic Casimir element $C$ of $L$ such
that $C_V \neq 0$\,. Since $V$ is simple, it follows that $C_V$ is invertible.
Thus Prop.~2.2 implies that $H^n(L,V) = \Ze$ if $n \neq 0$\,, and the same is
true for $n = 0$ since $0$ is the sole $L$--invariant element of $V$.

Thus we may now assume that $V$ is the trivial $L$--module $\KK$ and that
$n \ge 0$\,. According to Prop.~2.1.a every cohomology class contains an
$L$--invariant element. Conversely, it follows from Eq.~\reff{delalt} that
every $L$--invariant cochain (with values in $\KK$\,) is a cocycle.
Moreover, if $g$ and $g'$ are cohomologous elements of $Z^n(L,\KK)$ and if $g$
and $g'$ are $L$--invariant, then Prop.~2.1.b tells us that $g-g'=\de^{n-1}b$
with an $L$--invariant cochain $b \in C^{n-1}(L,\KK)$. As we just have seen,
this implies that $\de^{n-1}b = 0$\,, i.e., that $g = g'$.

Summarizing, we have proved that {\it $H^n(L,\KK)$ with $n \ge 1$ can be
identified with the space of all $L$--invariant skew symmetric $n$--linear
forms on $L$.}

The cohomology groups $H^n(L,\KK)$ with $n \in \{0,1,2\}$ can easily be
determined. Obviously, we have
 \be  H^0(L,\KK) = \KK\;\;{}  \ee
(see Eq.~\reff{hzero}). Let us next recall that the commutator algebra
$[L,L]$ of $L$ is equal to $L$\,. In view of Eq.~\reff{del1}, this implies
that $Z^1(L,\KK) = \Ze$ and hence
 \be  H^1(L,\KK) = \Ze \,.  \ee
Moreover, it also follows from $[L,L] = L$ that every $L$--invariant bilinear
form on $L$ is symmetric. This shows that
 \be  H^2(L,\KK) = \Ze \,.  \ee
On the other hand, it is easy to see that
 \be  H^3(L,\KK) \neq \Ze \,.  \ee
In fact, if $\phi$ is the Killing form of $L$\,, the assignment
 \be  (A,B,C) \lra \phi([A,B],C) \label{tri} \ee
(with $A,B,C \in L$) defines a non-zero, $L$--invariant, skew symmetric,
trilinear form on $L$\,.
\vspace{1ex}

\noindent
{\bf 2. The Lie superalgebras $osp(1|2n)$}\\[0.8ex]
For (finite--dimensional) simple Lie superalgebras the situation is much
more complicated, even for the basic classical ones. However, there is one
particular class of algebras (already mentioned by Kac in Sec.~5.5.3 of
Ref.~\cite{Ka}), namely, the $osp(1|2n)$ algebras, for which the discussion
of the preceding example goes through almost verbatim. All we have to do is
to interpret it in the super sense and to recall the bijection \reff{hshift}.
The pertinent information on the $osp(1|2n)$ algebras is contained in
Ref.~\cite{Dj}. As already mentioned, the cohomology groups $H^n(L,\KK)$ of
the non--exceptional classical simple Lie superalgebras have been described
in Ref.~\cite{FL} (see also Chap.~2, \S\,6 of Ref.~\cite{Fu} and
Ref.~\cite{Gr}; in particular, the exceptional algebras are considered in the
latter article).
\vspace{1ex}

\noindent
{\bf 3. Trivial coefficients}\\[0.8ex]
Let us return to the general assumptions (where $L$ is an arbitrary $\ve$
Lie algebra), and let us consider the case where $V = \KK$ is the trivial
graded $L$--module. As noted in Example 1, Eq.~\reff{delalt} implies that
every $L$--invariant cochain is a cocycle. We would like to mention a
classical construction of such cocycles (which, in the case of simple Lie
algebras, essentially does the whole job, see Ref.~\cite{GHV} for more
details).

Let $m \ge 0$ be an integer and let $\phi \in \Lg_{m+1}(L,\KK)$ be an
$L$--invariant $(m+1)$--linear form. We want to construct an
$L$--invariant cochain out of $\phi$\,. An obvious ansatz is to consider the
$\ve$--skew--symmetrization of the $2(m+1)$--linear form on $L$ defined by
 \be  (A_0,A_1,\ldots,A_{2m+1})
 \lra \phi(\la A_0,A_1\ra,\la A_2,A_3\ra,\ldots,\la A_{2m},A_{2m+1}\ra)
                                                        \label{skewsym} \ee
for all $A_i \in L$\,. However, the $L$--invariance of $\phi$ and the $\ve$
Jacobi identity imply, that the $\ve$--skew--symmetrization of the form
\reff{skewsym} is equal to zero. The well--known way out is to consider for
$m \ge 1$ the $\ve$--skew--symmetrization of the $(2m+1)$--linear form
$\psi \in \Lg_{2m+1}(L,\KK)$ defined by
 \be  \psi(A_0,A_1,\ldots,A_{2m})
  = \phi(\la A_0,A_1\ra,\la A_2,A_3\ra,\ldots,\la A_{2m-2},A_{2m-1}\ra,A_{2m})
                                                                         \ee
for all $A_i \in L$\,.

The situation simplifies a little if $\phi$ is $\ve$--symmetric (and
$L$--invariant). In that case it is sufficient to $\ve$--skew--symmetrize
$\psi$ in $A_1,A_2,\ldots,A_{2m-1}$\,, i.e., already the $(2m+1)$--linear
form on $L$ defined by
 \bea  (A_0,A_1,\ldots,A_{2m})
  &\! \lra \!& \!\!\sum_{\pi \in {\cal S}_{2m-1}}\!\!\mbox{sgn}(\pi)
                               \,\ve_{2m-1}(\pi;\al_1,\ldots,\al_{2m-1}) \\
    &      & \;\;\phi(\la A_0,A_{\pi(1)}\ra,\la A_{\pi(2)},A_{\pi(3)}\ra,
                   \ldots,\la A_{\pi(2m-2)},A_{\pi(2m-1)}\ra,A_{2m}) \nn \eea 
(for homogeneous elements $A_i \in L$) is $\ve$--skew--symmetric and, of
course, $L$--invariant. Note that the trilinear form \reff{tri} is the
simplest example of the foregoing construction.

The possibility to obtain $L$--invariant cochains from $L$--invariant
$\ve$--symmetric multilinear forms is welcome and interesting, since the
latter forms make their appearance in quite a different context. In fact,
if $L$ is finite--dimensional and admits a non--degenerate $L$--invariant
bilinear form which is homogeneous of degree zero, then there is a close
relationship between the Casimir elements of $L$ and the $\ve$--symmetric
$L$--invariant multilinear forms on $L$\,. For more details, we refer the
reader to Ref.~\cite{S.ce}. However, we note that we don't know of any
reason to hope that the situation might be as nice as for semi--simple Lie
algebras.
\vspace{1ex}  

\noindent
{\bf 4. Irreducible $gl(m|n)$--modules for which all Casimir operators \\
\hspace*{0.9em} are equal to zero}\\[0.8ex]
In this example we assume that the field $\KK$ is {\em algebraically closed.}
We want to determine all those finite--dimensional simple graded
$gl(m|n)$--modules, for which all Casimir operators without a constant term
are equal to zero. According to Prop.~2.2, among the finite--dimensional
simple graded modules, only these might have a non--trivial cohomology. In
view of the detailed knowledge of the Casimir elements of $gl(m|n)$ and of
their eigenvalues in a highest weight module \cite{S.je,S.ce,S.cec,S.ev} (see
also \cite{Ka.ca,Ka.bo,Ka.pr}), it is not difficult to solve this problem.

We are going to use the notation and results of the cited references. To
begin with, let us recall some of them. Let $m$ and $n$ be any strictly
positive integers. The Lie superalgebra $gl(m|n)$ is defined through its
elementary representation in $\KK^{m+n}$. If $(e_i)_{1 \le i \le m+n}$ is
the canonical basis of $\KK^{m+n}$, we define the gradation of $\KK^{m+n}$
by demanding that the elements $e_1,\ldots,e_m$ are even and the elements
$e_{m+1},\ldots,e_{m+n}$ are odd. Correspondingly, we define
 \be  \si_i \,=\, \Bigg\{ \begin{array}{r@{\;\;\;\mbox{if}\;\;\;}l}
                                    1 & 1 \le i \le m \\
                                   -1 & m+1 \le i \le m+n \,.
                          \end{array}                               \ee
In Ref.~\cite{S.ev} certain generators $X_{ij}$ of the Lie superalgebra
under consideration have been introduced. For $gl(m|n)$, these are just the
usual basic $(m+n) \ti (m+n)$ matrices $E_{ij}$\,; $1 \le i,j \le m+n$\,.

Now let $s \ge 1$ be an integer and let $C_s$ be the element of the universal
enveloping algebra $U(gl(m|n))$ defined by
 \be  C_s = \sum_{i_1,\ldots,\,i_s}\si_{i_1}\si_{i_2}\ldots\si_{i_{s-1}}
                      \,X_{i_s i_1}X_{i_1 i_2}\ldots X_{i_{s-1}i_s} \,. \ee
Then $C_s$ is a Casimir element of $gl(m|n)$, and (as shown in
Ref.~\cite{S.cec}) the super center $Z(gl(m|n))$ of $U(gl(m|n))$ is
generated (as an algebra) by the unit element and these elements $C_s$\,.
Stated differently, every Casimir element can be written as a polynomial in
the $C_s$\,; $s \ge 1$\,. (By the way, none of the $C_s$ may be dropped. In
particular, the algebra $Z(gl(m|n))$ is not finitely generated.)

The elements $X_{ii}$\,; $1 \le i \le m+n$ span the standard Cartan
subalgebra $H$ of $gl(m|n)$. Let $\ve_i$\,; $1 \le i \le m+n$ be the linear
form on $H$ which satisfies
 \be  \ve_i(X_{jj}) = \de_{ij}  \ee
for $1 \le j \le m+n$\,. We choose the usual system of positive roots of
$gl(m|n)$ with respect to $H$, namely
 \be  \De^+ = \{\ve_i - \ve_j \,|\, 1 \le i < j \le m+n\} \,.  \ee

Now let $\Lam \in H^{\ast}$ be any linear form on $H$ and let $V$ be a graded
highest weight $gl(m|n)$--module with highest weight $\Lam$\,. If
$C \in Z(gl(m|n))$ is any Casimir element, the corresponding Casimir operator
$C_V$ is a scalar multiple of $id_V$\,,
 \be  C_V = \chi_{\Lam}(C)\,id_V \,,  \ee
where $\chi_{\Lam}$ is the so--called central character associated to
$\Lam$\,. To solve our problem, we need the eigenvalues $\chi_{\Lam}(C)$.
According to the foregoing remarks, it is sufficient to know the
$\chi_{\Lam}(C_s)$\,; $s \ge 1$\,. These have been investigated in
Ref.~\cite{S.ev}. To describe the pertinent result, we introduce some more
notation.

Let $2\rho$ be the sum of the even positive roots minus the sum of the odd
positive roots, i.e.,
 \be  \rho = \frac{1}{2} \sum_{i < j}\si_i \si_j (\ve_i - \ve_j) \,.  \ee
We set 
 \be  \lam = \Lam + \rho  \ee
and define the numbers $r_i$ and $\ell_i$ through
 \bea  r_i \eqq \si_i \,\rho(X_{ii}) \label{defr} \\
    \ell_i \eqq \si_i \,\lam(X_{ii}) \,. \label{deflam} \eea
If $C \in Z(gl(m|n))$ is fixed, we have
 \be  \chi_{\Lam}(C) = P_C(\ell_1,\ldots,\ell_{m+n}) \,,  \ee
where $P_C$ is a polynomial in $m+n$ indeterminates $Y_1,\ldots,Y_{m+n}$\,.
If $C$ runs through all of $Z(gl(m|n))$, the polynomials $P_C$ form an
algebra $T(m|n)$ (the algebra of supersymmetric polynomials in super
dimension $(m|n)$ \cite{S.je}; for a recent reference see Ref.~\cite{Mo}),
and the map of $Z(gl(m|n))$ onto $T(m|n)$ defined by $C \rar P_C$
is an algebra isomorphism (the generalized Harish--Chandra isomorphism
\cite{S.ce}). Finally, the algebra $T(m|n)$ is generated by the unit element
and the polynomials
 \be  Q_s = \sum_{i=1}^{m+n}\si_i(Y_i^s - r_i^s)  \ee
with $s \ge 1$ (see the end of Sec.~2 in Ref.~\cite{S.ev}, in particular, the
lines below Eq.~(2.36) and Eq.~(2.37)). Note that the constant term of a
Casimir element $C$ is equal to $\chi_0(C)$. Consequently, this term is
equal to zero if and only if $P_C$ can be written as a polynomial in the
$Q_s$ without a constant term.

Next we recall that every finite--dimensional simple graded $gl(m|n)$--module
is a highest weight module (this is the sole place where we use that the base
field is algebraically closed). Summarizing, we conclude that our problem can
be reformulated as follows: \\
Determine all linear forms $\Lam \in H^{\ast}$ which correspond to
finite--dimensional simple graded $gl(m|n)$--modules and satisfy
 \be  Q_s(\ell_1,\ldots,\ell_{m+n}) = 0 \label{qvan}  \ee
for all $s \ge 1$\,.

Let us set
 \be  L_i = \Lam(X_{ii})  \ee
for $1 \le i \le m+n$\,. (Note that in contrast to Eqs.~\reff{defr},
\reff{deflam} we have not included a factor of $\si_i$ on the right hand
side.) Then the first condition is satisfied if and only if
 \be  L_i - L_{i+1} \in \NN \label{findim}  \ee
for $1 \le i \le m-1$ and for $m+1 \le i \le m+n-1$\,.

The validity of the relations \reff{qvan} is easily seen to be equivalent
to the following condition: \\
{\it Up to the ordering, the numbers
 \be  \ell_1,\ell_2,\ldots,\ell_m,r_{m+1},r_{m+2},\ldots,r_{m+n}  \ee
coincide with the numbers
 \be  r_1,r_2,\ldots,r_m,\ell_{m+1},\ell_{m+2},\ldots,\ell_{m+n} \,.  \ee }

The subsequent discussion depends on whether $m = n$\,, $m > n$\,, or
$m < n$\,. Using Eq.~\reff{findim} and some special properties of the
numbers $r_i$\,, we first conclude that the conditions \reff{qvan} are
satisfied if and only if the following holds true. \\
If $m = n$\,, we have
 \be \ell_{m+i}\: = \:\ell_{m+1-i}    \hspace{1em}\mbox{for}\;\;\;
                                             1 \le i \le m \,,\hst{1.65em} \ee
if $m > n$\,, there exists a number $k \in \{0,1,\ldots,n\}$ such that
\bea \ell_{m+i} \eqq \ell_{m+1-i} \hspace{1em}\mbox{for}\;\;\;
                                                    1 \le i \le k \\
        r_{j-n} \eqq \ell_{j-k}   \hspace{2.11em}\mbox{for}\;\;\;
                                             n+1 \le j \le m \hst{0.27em}\\
     \ell_{m+i} \eqq \ell_{n+1-i} \hspace{1.2em}\mbox{for}\;\;\;
                                                    k+1 \le i \le n \,, \eea
and if $m < n$\,, there exists a number $h \in \{0,1,\ldots,m\}$ such that
\bea \ell_{m+i} \eqq \ell_{m+1-i} \hspace{1em}\mbox{for}\;\;\;
                                                    1 \le i \le h \\
     \ell_{h+j} \eqq r_{m+j}      \hspace{1.82em}\mbox{for}\;\;\;
                                                    m+1 \le j \le n \\
     \ell_{n+i} \eqq \ell_{m+1-i} \hspace{1em}\mbox{for}\;\;\;
                                                    h+1 \le i \le m \,. \eea
Expressed in terms of the coefficients $L_i$\,, our final result reads as
follows.

Suppose first that $m = n$\,. Choose the numbers $L_1,L_2,\ldots,L_m$ in
accordance with Eq.~\reff{findim} and define
 \be  L_{m+i} = -L_{m+1-i} \hspace{1em}\mbox{for}\;\;\;
                                                    1 \le i \le m \,.  \ee
Suppose next that $m > n$\,. Let $k$ be any element of $\{0,1,\ldots,n\}$.
Choose the numbers $L_1,L_2,\ldots,L_m$ in accordance with Eq.~\reff{findim}
and such that
 \be  L_{n+1-k} = L_{n+2-k} = \ldots = L_{m-k} = n - k \,.  \ee
Define
\bea
 L_{m+i} \eqq -L_{m+1-i}       \hspace{5.56em}\mbox{for}\;\;\; 1 \le i \le k \\
 L_{m+i} \eqq -L_{n+1-i} - (m-n) \hspace{1em}\mbox{for}\;\;\; k+1 \le i \le n
                                                                        \eea
(for $k = 0$ resp. $k = n$\,, the first resp. second of these equations
should be dropped). \\
Finally, suppose that $m < n$\,. Let $h$ be any element of $\{0,1,\ldots,m\}$.
Choose the numbers $L_{m+1},L_{m+2},\ldots,L_{m+n}$ in accordance with
Eq.~\reff{findim} and such that
 \be  L_{m+1+h} = L_{m+2+h} = \ldots = L_{n+h} = -(m - h) \,.  \ee
Define
\bea
 L_i \eqq -L_{m+n+1-i} + (n-m) \hspace{1em}\mbox{for}\;\;\; 1 \le i \le m-h \\
 L_i \eqq -L_{2m+1-i}          \hspace{6.36em}\mbox{for}\;\;\;
                                                      m-h+1 \le i \le m \eea
(for $h = m$ resp. $h = 0$\,, the first resp. second of these equations
should be dropped). \\
Then, in all three cases, our conditions are satisfied. Conversely, every
solution to our conditions is of this form.

We note that the relation \reff{qvan} with $s = 1$ is equivalent to
 \be  \sum_{i=1}^{m+n}L_i = 0 \,, \label{strvan}  \ee
and it is easily checked that the numbers $L_i$ specified above satisfy
this condition.

\noindent
{\it Remark 3.1.} Needless to say the modules we have found are atypical
\cite{Ka.bo}. Actually, Eqs.~(3.25)\,--\,(3.31) show that they are maximally
atypical in the sense that $\Lam$ satisfies the maximal number (equal to
$\min(m,n)$) of atypicality conditions.

The foregoing results can immediately be extended to the special linear Lie
superalgebras $sl(m|n)$, provided that $m \neq n$\,. In fact, all we have to
do is to interpret them appropriately, as explained in Ref.~\cite{S.ev}. The
elements $X_{ij}$ are now given by
 \be  X_{ij} = E_{ij} - \frac{1}{d}\,\si_i\,\de_{ij}\,I \,,  \ee
where we have set $d = m-n$\,, and where $I$ is the unit matrix. Furthermore,
the linear forms $\ve_i$ are given by
 \be  \ve_i(X_{jj}) = \de_{ij} - \frac{1}{d}\,\si_j \,.  \ee
This implies that
 \be  \sum_{i=1}^{m+n}X_{ii} = 0  \label{slgen}\ee
and
 \be  \sum_{i=1}^{m+n}\si_i \,\ve_i = 0 \,.  \ee
Apart from these modifications, everything goes through, and the highest
weights of the modules we are looking for are given by the same formulae as
above. (We have already noted that the relation \reff{strvan}, which should
be satisfied because of Eq.~\reff{slgen}, indeed holds true.)
\vspace{5ex}

%
\noindent
{\Large \bf
IV. Central extensions of $\ve$ Lie algebras \\[-0.5ex]}
\setcounter{section}{4} \setcounter{equation}{0}

\noindent
We start with the basic definitions. Let $L$ be an $\ve$ Lie algebra. A
{\it central extension} of $L$ by an $\ve$ Lie algebra $H$ is an exact
sequence of $\ve$ Lie algebras,
 \be  \Ze \lra H \ario E \arpi L \lra \Ze \, \label{cen}  \ee
such that $\io(H)$ is in the $\ve$ center of $E$\,. Thus $\io$ and $\pi$
are homomorphisms of $\ve$ Lie algebras (which implies that they are
homogeneous of degree zero), $\io$ is injective, $\pi$ is surjective, and
the kernel of $\pi$ is equal to the image of $\io$\,. In particular, $H$ is
Abelian. Thus $H$ is nothing but a graded vector space, endowed with the
trivial structure of an $\ve$ Lie algebra. If
 \be  \Ze \lra H' \ariopr E' \arpipr L \lra \Ze \label{cenpr} \ee
is a second central extension of $L$\,, then by a {\it morphism} from the
central extension \reff{cen} to the central extension \reff{cenpr}, we
mean a pair $(\vp,\psi)$ of $\ve$ Lie algebra homomorphisms,
 \be  \vp : E \lra E' \;\;\;,\;\;\; \psi : H \lra H'  \ee
such that
 \be \pi = \pi' \comp \vp \;\;\;,\;\;\; \vp \comp \io = \io' \comp \psi \,, \ee
i.e., such that the corresponding diagram is commutative. (Obviously, $\psi$
is uniquely determined by $\vp$\,.) The morphism $(\vp,\psi)$ is said to be
an {\it isomorphism} of the central extensions if $\vp$ and $\psi$ are
bijective. Actually, if one of the maps $\vp$ and $\psi$ is bijective, then
so is the other.

Consider next the case where the second central extension has the form
 \be  \Ze \lra H \ariopr E' \arpipr L \lra \Ze \label{cenprsp} \,.  \ee
Then the extensions \reff{cen} and \reff{cenprsp} are said to be
{\it equivalent} if there exists a morphism $(\vp,\psi)$ from the extension
\reff{cen} to the extension \reff{cenprsp} such that $\psi = id_H$ (which
implies that $\vp$ is bijective).

As in the Lie algebra case, the central extensions of $L$ by $H$ can be
classified by means of the cohomology groups $H^2(L,H)$, where $H$ is
regarded as a graded vector space, endowed with the trivial $L$--module
structure. To show this, we choose a homogeneous {\em section} $\si$
(in the category of vector spaces) of
$\pi$\,, i.e., a homogeneous linear map
 \be  \si : L \lra E  \ee
(necessarily of degree zero) such that
 \be  \pi \comp \si = id_L \,.  \ee
Then we have
 \be  \pi(\la\si(A),\si(B)\ra - \si(\la A,B\ra)) = 0  \ee
for all $A,B \in L$\,. Consequently, there exists a unique map
 \be  g : L \ti L \lra H  \ee
such that
 \be  \io(g(A,B)) = \la\si(A),\si(B)\ra - \si(\la A,B\ra) \label{defcocy} \ee
for all $A,B \in L$\,, and it is easy to see that $g$ is a homogeneous
2--cocycle of degree zero, i.e., $g \in Z^2(L,H)_0$\,.

Consider in addition the central extension \reff{cenprsp}, choose a
homogeneous section $\si'$ of $\pi'$, and construct the corresponding
2--cocycle $g'$. Then the central extensions \reff{cen} and \reff{cenprsp}
are equivalent if and only if $g$ and $g'$ are cohomologous.

Conversely, let $g$ be any element of $Z^2(L,H)_0$\,. Make the graded vector
space $L \ti H$ into an $\ve$ Lie algebra by defining
 \be  \la(A,X),(B,Y)\ra = (\la A,B\ra,g(A,B))  \ee
for all $A,B \in L$ and all $X,Y \in H$. This $\ve$ Lie algebra will be
denoted by $L(g)$. Let
 \be  \si : L \lra L \ti H \;\;\;,\;\;\; \io : H \lra L \ti H  \ee
be the canonical injections and let
 \be  \pi : L \ti H \lra L  \ee
be the canonical projection. Then
 \be  \Ze \lra H \ario L(g) \arpi L \lra \Ze \label{ceng} \ee
is a central extension of $L$ by $H$, $\si$ is a homogeneous section of
$\pi$\,, and the corresponding 2--cocycle of $L$ with values in $H$ is just
the 2--cocycle $g$ we started with.

Summarizing, we conclude that there is a bijection between the set of
equivalence classes of central extensions of $L$ by $H$ and the cohomology
group $H^2(L,H)$. (More generally, an analogous result holds for all those
extensions \reff{cen} of $L$ for which $H$ is only assumed to be Abelian;
for more details, see Ref.~\cite{SSL}.)

In the Lie algebra case, all this is well--known, and it is not at all
surprising that this material can be generalized to the present graded
setting. We have included these results in order to fix our notation and as
a preparation for the investigations to follow. In fact, we now are ready to
generalize Garland's theory of universal central extensions of Lie algebras
\cite{Ga}. As we are going to see, most of Garland's statements and proofs
can be adopted almost verbatim. In addition, we comment on the
{\it homological}\/ background of these results, and we show how universal
central extensions can be constructed.

Let us agree to call an $\ve$ Lie algebra $L$ {\it perfect}\/ if it is equal
to its own commutator algebra, i.e., if $\la L,L\ra = L$\,. Using this
notion, we call the central extension \reff{cen} a {\it covering}\/ of $L$
if $E$ is perfect. In this case, we also call $E$ or $\pi$ or the pair
$(E,\pi)$ a covering of $L$\,, and we shall say that $E$ (or $\pi$ or
$(E,\pi)$) covers $L$\,. Obviously, if an $\ve$ Lie algebra admits a covering,
then it is perfect.

\noindent
{\bf Lemma 4.1}\\
{\it If the central extension \reff{cen} is a covering of $L$\,, then there
exists at most one morphism from \reff{cen} to any second central extension
of $L$\,. }

\noindent
Proof \\
Let $(\vp,\psi)$ and $(\vp',\psi')$ be two morphisms from the central
extension \reff{cen} to the central extension \reff{cenpr}. Obviously, we
have
 \be  (\vp - \vp')(\la A,B\ra)
    = \la\vp(A) - \vp'(A),\vp(B)\ra + \la\vp'(A),\vp(B) - \vp'(B)\ra
                                                         \label{diff}  \ee 
for all $A,B \in E$\,. On the other hand, we have
 \be  \pi'(\vp(A) - \vp'(A)) = \pi(A) - \pi(A) = 0  \ee
and hence
 \be  \vp(A) - \vp'(A) \in \io'(H')  \ee
for all $A \in E$\,. Consequently, the right hand side of Eq.~\reff{diff}
vanishes, which implies that $\vp = \vp'$, and hence also $\psi = \psi'$.

\noindent
{\bf Definition 4.1}\\
{\it We say that a covering of an $\ve$ Lie algebra $L$ is {\em universal,}
if for every central extension of $L$ there is a unique morphism (in the
sense of central extensions) from the covering to the central extension.}

We note that the uniqueness requirement in Def.~4.1 automatically follows
from Lemma 4.1. Obviously, the definition and Lemma 4.1 imply the following
proposition.

\noindent
{\bf Proposition 4.2}\\
{\it Any two universal coverings of an $\ve$ Lie algebra $L$ are isomorphic
as central extensions, and the two mutually inverse isomorphisms are unique.}
 
We have already noted that if an $\ve$ Lie algebra admits a covering, then
it is perfect. Conversely, we have:

\noindent
{\bf Proposition 4.3}\\
{\it If the $\ve$ Lie algebra $L$ is perfect, then $L$ has a universal
covering.}

\noindent
Proof \\
Let $\Wept L$ be the $\ve$ exterior algebra of the graded vector space $L$
(see App.~A for more information), let $\wepe$ denote the multiplication in
$\Wept L$\,, and let $W' = \Weptt L$ be the subspace of $\Wept L$ of
$\ZZ$--degree 2. Consider the subspace $I$ of $W'$ which is generated by all
elements of the form
 \be  - \la A,B\ra \wepf C + \ve(\bet,\ga)\,\la A,C\ra \wepf B
                                             + A \wepf \la B,C\ra \,,  \ee
with homogeneous elements $A,B,C \in L$\,. Obviously, $I$ is a graded
subspace of $W'$. We set $W = W'/I$ and denote the canonical image of
$A \wepf B$ in $W$ by $f(A,B)$, for all $A,B \in L$\,. By definition,
 \be  f : L \ti L \lra W  \ee
is a 2--cocycle which is homogeneous of degree zero. Thus we have the central
extension
 \be  \Ze \lra W \ariob L(f) \arpib L \lra \Ze \,. \label{cenf}  \ee

Consider an arbitrary central extension of $L$ by an Abelian $\ve$ Lie
algebra $H$. We know that this extension is equivalent to the extension
\reff{ceng} with a suitable 2--cocycle $g \in Z^2(L,H)_0$\,. Since $g$ is
a 2--cocycle, there exists a unique linear map
 \be  \psi' : W \lra H  \ee
such that
 \be  \psi'(f(A,B)) = g(A,B)  \ee
for all $A,B \in L$\,, and $\psi'$ is homogeneous of degree zero. We define
a map
 \be  \vp' : L(f) = L \ti W \lra L \ti H = L(g)  \ee
by
 \be  \vp'(A,X) = (A,\psi'(X))  \ee
for all $A \in L$ and $X \in W$. Then $(\vp',\psi')$ is a morphism from
the central extension \reff{cenf} to the central extension \reff{ceng}.

Now let $\Lh = \la L(f),L(f)\ra$ be the commutator algebra of $L(f)$.
Since $L$ is perfect, we have $\Lh + W = L(f)$ and hence
 \be  \bar{\pi}(\Lh) = L  \ee
 \be  \Lh = \la \Lh+W,\Lh+W \ra = \la \Lh,\Lh \ra \,. \ee
Thus $\Lh$ is perfect. If $\Hh$ is the preimage of $\Lh$ under
$\bar{\io}$\,, the central extension \reff{cenf} induces a central extension
 \be \Ze \lra \Hh \arioh \Lh \arpih L \lra \Ze \,, \label{cenuni} \ee 
which is a covering of $L$\,. Furthermore, if $\vp$ is the restriction of
$\vp'$ to $\Lh$ and if $\psi$ is the restriction of $\psi'$ to $\Hh$,
the pair $(\vp,\psi)$ is a morphism from the central extension \reff{cenuni}
to the central extension \reff{ceng}. This shows that \reff{cenuni} is a
universal covering of $L$\,, and thus the proposition is proved.

In the present work we do not consider the homology of $\ve$ Lie algebras
(it is fairly obvious how this should be introduced). Nevertheless, it would
be rather unsatisfactory not to mention the homological interpretation of the
preceding proof.

Consider the linear map
 \be  \pa_2 : \Wepft L \lra L  \ee
given by
 \be  \pa_2(A \wepf B) = -\la A,B\ra \,,  \ee
and the linear map
 \be  \pa_3 : \Wepfd L \lra \Wepft L  \ee
specified by
 \be  \pa_3(A \wepf B \wepf C)
 = -\la A,B\ra \wepf C + \ve(\bet,\ga)\,\la A,C\ra \wepf B + A \wepf \la B,C\ra
                                                                     \,,  \ee
for all homogeneous elements $A,B,C \in L$\,. The $\ve$ Jacobi identity is
equivalent to
 \be  \pa_2 \comp \pa_3 = 0 \,. \label{papa}  \ee
Let $Z_2(L)$ be the kernel of $\pa_2$ and let $B_2(L)$ be the image of
$\pa_3$\,. Then Eq.~\reff{papa} implies that $B_2(L) \subset Z_2(L)$, and
a closer look at the definition of $\Hh$ in the proof of Prop.~4.3 shows
that
 \be  \Hh = Z_2(L)/B_2(L) = H_2(L)  \ee
(where the last equation is a definition). By definition, $C^n(L,\KK)$ is the
graded dual of $\Weptn L$ (see the remark below Eq.~\reff{cochz}), and the
coboundary map
 \be  \de^n : C^n(L,\KK) \lra C^{n+1}(L,\KK)  \ee
is the graded transpose of the boundary map
 \be  \pa_{n+1} : \Wepfnpe L \lra \Wepfn L  \ee
(we are considering the cases $n = 1,2$ only, but, in fact, this is true for
all $n$). Consequently, the graded dual $\Hhastgr$ of $\Hh$ is
canonically isomorphic to $H^2(L,\KK)$.

More explicitly, the canonical isomorphism
$\phi: H^2(L,\KK) \rar \Hhastgr$ can be described as follows. Let
$g \in Z^2(L,\KK)_{\ga}$ be a cocycle which is homogeneous of degree $\ga$\,.
We have to give the image of the cohomology class $g + B^2(L,\KK)$ in
$\Hhastgr$. The cocycle $g$ can be regarded as an element
$g_0 \in Z^2(L,\KK^{\ga})_0$ (recall Eq.~\reff{grshift}). Let
 \be  \Ze \lra \KK^{\ga} \lra E \lra L \lra \Ze \label{cengz} \ee
be a central extension of $L$ which belongs to the equivalence class
corresponding to $g_0 + B^2(L,\KK^{\ga})$, and let $(\vp_g,\psi_g)$ be the
morphism from the universal covering \reff{cenuni} to \reff{cengz}. Then
$\psi_g$ can be regarded as a linear form $\Hh \rar \KK$ which is
homogeneous of degree $\ga$\,. This form depends only on the cohomology
class $g + B^2(L,\KK)$ and is the image we are looking for.

The inverse of the isomorphism $\phi$ can be described even more simply. Let
$\hat{g} \in Z^2(L,\Hh)_0$ be a cocycle such that the cohomology class
$\hat{g} + B^2(L,\Hh)$ corresponds to the equivalence class of the
universal covering \reff{cenuni}. If $\lam$ is any element of
$\Hhastgr$, then $\lam \comp \hat{g}$ is an element of $Z^2(L,\KK)$
and we have
 \be  \phi^{-1}(\lam) = \lam \comp \hat{g} + B^2(L,\KK) \,.  \ee

If $L$ is a perfect $\ve$ Lie algebra for which $H^2(L,\KK)$ is
finite--dimensional (of dimension $s$\,, say), the preceding discussion leads
to the following construction of a universal covering of $L$\,. Let
$(g_r)_{1 \le r \le s}$ be a family of homogeneous 2--cochains
$g_r \in Z^2(L,\KK)$, such that the cohomology classes $g_r + B^2(L,\KK)$
form a basis of $H^2(L,\KK)$, and let $\ga_r$ denote the degree of $g_r$\,.
Choose a graded vector space $H$ such that $H$ has a homogeneous basis
$(e_r)_{1 \le r \le s}$\,, with $e_r$ homogeneous of degree $-\ga_r$\,.
Define the bilinear map
 \be  g : L \ti L \lra H  \ee
by
 \be  g(A,B) = \sum_{r=1}^{s}g_r(A,B)\,e_r \label{univcocy} \ee
for all $A,B \in L$\,. Then $g \in Z^2(L,H)_0$\,, and the corresponding
central extension \reff{ceng} is a universal covering of $L$\,.

\noindent
{\it Remark 4.1.} Let
 \be \Ze \lra \Hh \arioh \Lh \arpih L \lra \Ze \label{covuni} \ee
be a universal covering of the $\ve$ Lie algebra $L$\,. Then $\Lh$ is a
universal covering of every $\ve$ Lie algebra $E$ which covers $L$\,. In
fact, if
 \be  \Ze \lra H \lra E \lra L \lra \Ze \label{cov} \ee
is a covering of $L$ and if $(\vp,\psi)$ is a morphism from \reff{covuni}
to \reff{cov}, then it is easy to see that
 \be  \Ze \lra \mbox{ker}(\vp) \lra \Lh \arvp E \lra \Ze  \ee
is a covering of $E$\,. One can then show directly that this covering is
universal, but the simplest way to prove this is to invoke Thm.~4.1.

\noindent
{\it Remark 4.2.} The universal coverings have a property which, apparently,
is more general than the one required in their definition. Let
 \be  \Ze \lra H \ario E \arpi L \lra \Ze  \ee
be a universal covering of an $\ve$ Lie algebra $L$ and let
 \be  \Ze \lra H' \ariopr E' \arpipr L' \lra \Ze  \ee
be a central extension of an $\ve$ Lie algebra $L'$. If $\vr : L \rar L'$
is a homomorphism (of $\ve$ Lie algebras), there exists a unique pair
$(\vp,\psi)$ of homomorphisms
 \be  \vp : E \lra E' \;\;\;,\;\;\; \psi : H \lra H'  \ee
such that
 \be  \vr \comp \pi = \pi' \comp \vp \;\;\;,\;\;\;
                                  \vp \comp \io = \io' \comp \psi \,,  \ee
i.e., such that the corresponding diagram is commutative.

\noindent
{\it Remark 4.3.} One might think to generalize the definition of a universal
covering to all (not necessarily perfect) $\ve$ Lie algebras, as follows: A
central extension \reff{cen} of an $\ve$ Lie algebra $L$ is said to be
universal if, for any central extension \reff{cenpr} of $L$\,, there exists
a unique morphism from \reff{cen} to \reff{cenpr}. Actually, it can be shown
that if \reff{cen} is universal in this sense, then $L$ and $E$ are perfect,
and nothing beyond the universal coverings has been obtained.

\noindent
{\bf Definition 4.2}\\
{\it a) An $\ve$ Lie algebra $L$ is said to be {\em simply connected,} if for
every central extension
 \be  \Ze \lra H \ario E \arpi L \lra \Ze  \ee
of $L$\,, there is a unique homomorphism of $\ve$ Lie algebras
$\si : L \rar E$\,, such that $\pi \comp \si = id_L$\,. \\
b) A covering
 \be  \Ze \lra H \lra E \lra L \lra \Ze  \ee
of an $\ve$ Lie algebra $L$ is said to be simply connected, in case $E$ is
simply connected.}

Every simply connected $\ve$ Lie algebra $L$ is perfect. In fact, choose
$H = L/\la L,L\ra$ and consider the trivial central extension
 \be  \Ze \lra H \lra L \ti H \lra L \lra \Ze  \ee
of $L$\,, where $H \rar L \ti H$ is the canonical injection and
$L \ti H \rar L$ is the canonical projection. Of course, the canonical
injection $L \rar L \ti H$ is a homomorphism $\si$ of the type described
in Def.~4.2.a. But if $\vp : L \rar L/\la L,L\ra$ is the canonical
homomorphism, the map
 \be  \si' : L \lra L \ti H  \ee
defined by
 \be  \si'(A) = (A,\vp(A))  \ee
for all $A \in L$ is likewise. Consequently, the uniqueness of $\si$ implies
our claim.

\noindent
{\bf Theorem 4.1}\\
{\it A covering of an $\ve$ Lie algebra $L$ is universal if and only if it is
simply connected.}

\noindent
Proof \\
Let
 \be  \Ze \lra \Hh \arioh \Lh \arpih L \lra \Ze \label{covthm} \ee
be a simply connected covering of $L$\,. We have to show that there exists
a (necessarily unique) morphism from this covering to an arbitrary central
extension \reff{cen} of $L$\,.

Define
 \be  E' = \{(A,B) \in E \ti \Lh \,|\, \pi(A) = \pih(B) \} \,.  \ee
Obviously, $E'$ is a graded subalgebra of the $\ve$ Lie algebra
$E \ti \Lh$\,. Let $\vr : E' \rar \Lh$ be the projection onto the second
factor of $E \ti \Lh$\,. Then $\vr$ is a surjective homomorphism of $\ve$
Lie algebras, and the kernel of $\vr$ is equal to
 \be  \mbox{ker}(\vr) = \{(\io(X),0)\,|\, X \in H\} \,,  \ee
which is contained in the center of $E'$. Since $\Lh$ is simply connected,
there exists a unique homomorphism of $\ve$ Lie algebras $\si : \Lh \rar E'$,
such that $\vr \comp \si = id_{\Lh}$\,. Let $\lam : \Lh \rar E$ denote the
composition of $\si$ with the projection of $E'$ onto the first factor of
$E \ti \Lh$\,. The definitions imply that
 \be  \si(B) = (\lam(B),B) \in E'  \ee
for all $B \in \Lh$ and hence (by definition of $E'$) that
$\pi \comp \lam = \pih$\,. In particular, we have
$\pi \comp \lam \comp \hat{\iota} = 0$\,, which shows that
 \be  (\lam \comp \hat{\iota})(\Hh) \subset \io(H) \,.  \ee
Consequently, there exists a unique homomorphism of $\ve$ Lie algebras
$\mu : \Hh \rar H$ such that $\lam \comp \hat{\iota} = \io \comp \mu$\,.
Then $(\lam,\mu)$ is a morphism from the covering \reff{covthm} to the
central extension \reff{cen}.

Conversely, let us assume that \reff{covthm} is a universal covering of $L$\,.
Consider a central extension
 \be  \Ze \lra H \ario E \arpi \Lh \lra \Ze \label{centhm} \ee
We have to show that there exists a unique homomorphism of $\ve$ Lie algebras
$\si : \Lh \rar E$ such that $\pi \comp \si = id_{\Lh}$\,.

If $\si$ exists, the pair $(\si,0)$ is a morphism from the covering
 \be  \Ze \lra \Ze \lra \Lh \arid \Lh \lra \Ze  \ee
to the central extension \reff{centhm}. Since $\Lh$ is perfect, Lemma 4.1
implies that $\si$ is uniquely determined. Hence it remains to prove that
$\si$ exists.

Let $\Et = \la E,E\ra$ be the commutator algebra of $E$\,. Since $\Lh$ is
perfect, we have $\pi(\Et) = \Lh$ and hence $E = \Et + \io(H)$. This in
turn implies that
 \be  \Et = \la \Et + \io(H),\Et + \io(H)\ra = \la\Et,\Et\ra \,,  \ee
i.e., that $\Et$ is perfect. Let $\Ht$ be the preimage of $\Et$ under $\io$\,.
Then the central extension \reff{centhm} induces a central extension
 \be  \Ze \lra \Ht \lra \Et \lra \Lh \lra \Ze \,, \label{cenred} \ee
and the homomorphism $\si$ exists, if the analogous homomorphism exists for
\reff{cenred}. Thus we may assume that $E$ is perfect.

Let $\nu : E \rar L$ be the (surjective) homomorphism of $\ve$ Lie algebras
defined by $\nu = \pih \comp \pi$\,. We are going to show that the
kernel of $\nu$ is contained in the center of $E$\,. Thus suppose that $A$
is a homogeneous element of $E$ such that $\pih(\pi(A)) = 0$\,. Then
$\pi(A)$ is central in $\Lh$\,. Consequently, we have
 \be  \pi(\la A,B\ra) = \la\pi(A),\pi(B)\ra = 0  \ee
for all $B \in E$\,, i.e., $\la A,B\ra$ belongs to the kernel of $\pi$ and
hence to the center of $E$\,. If $B$ and $B'$ are any two homogeneous
elements of $E$\,, this implies that
 \be  \la A,\la B,B'\ra\ra = \la\la A,B\ra,B'\ra
                          + \ve(\al,\bet)\,\la B,\la A,B'\ra\ra = 0 \,, \ee
and since $E$ is perfect it follows that $A$ lies in the center of
$E$\,.

Thus
 \be  \Ze \lra \mbox{ker}(\nu) \lra E \arnu L \lra \Ze \label{cennu} \ee
is a central extension. Since we are assuming that \reff{covthm} is a
universal covering, there exists a morphism $(\vp,\psi)$ from \reff{covthm}
to \reff{cennu}. We are going to show that $\vp$ is the homomorphism $\si$
we are looking for. To do so we have to show that $\pi \comp \vp = id_{\Lh}$\,.
Let $\om : \Lh \rar \Lh$ be the linear map defined by
$\om = \pi \comp \vp - id_{\Lh}$\,. Since $(\vp,\psi)$ is a morphism of
central extensions, we have
 \be  \pih \comp \om = \pih \comp (\pi \comp \vp) - \pih
                      = \nu \comp \vp - \pih = 0 \,,  \ee
and hence $\om$ takes its values in the center of $\Lh$\,. This implies that
 \be  \om(\la A,B\ra) = \la\om(A),(\pi \comp \vp)(B)\ra
                                                  + \la A,\om(B)\ra = 0  \ee
for all $A,B \in \Lh$\,. Since $\Lh$ is perfect, it follows that $\om = 0$
and the theorem is proved.
\vspace{1ex}

\noindent
{\bf Example: The universal covering of $sl(n|n)/\KK \cd I_{2n}$ }\\[0.8ex]
For abbreviation, we write $L$ instead of $sl(n|n)/\KK \cd I_{2n}$ (here and
in the following, $I_{2n}$ denotes the $2n \ti 2n$ unit matrix). It is
well--known \cite{S} that $L$ inherits from $gl(n|n)$ a consistent
$\ZZ$--gradation such that
 \be  L = L_{-1} \oplus L_0 \oplus L_1 \,.  \ee
To begin with, we note that
 \be  \Ze \lra \KK \ario sl(n|n) \arpi L \lra \Ze \label{slcov}  \ee
is a central extension of $L$\,, where $\io$ is defined by
 \be  \io(c) = \frac{c}{2n}\,I_{2n}  \ee
for all $c \in \KK$\,, and where $\pi$ is the canonical homomorphism.
Obviously, there exists a unique linear map
 \be  \si : L \lra sl(n|n)  \ee
such that
 \be  \si(\pi(A)) = A - \frac{1}{2n}\,\mbox{Tr}(A)\,I_{2n}  \ee
for all $A \in sl(n|n)$, and we have $\pi \comp \si = id_L$\,. It follows that
the 2--cocycle $g \in Z^2(L,\KK)$ defined by Eq.~\reff{defcocy} is given by
 \be  g(\pi(A),\pi(B)) = \mbox{Tr}(\la A,B\ra)  \ee
for all $A,B \in sl(n|n)$. Obviously, $g$ is $L_0$--invariant (but, of course,
not $L$--invariant).

If $n = 1$\,, the Lie superalgebra $L$ is Abelian and hence the theory of
universal coverings does not apply (see Remark 4.3). Thus from now on we shall
assume that $n \ge 2$\,. Then $L$ is simple and hence perfect. Our main task is
to determine $H^2(L,\KK)$. Since $L_0 \simeq sl(n) \ti sl(n)$ is semi--simple,
Prop.~2.1 implies that every cohomology class contains an $L_0$--invariant
cocycle. On the other hand, there does not exist a non--zero $L_0$--invariant
linear form on $L$\,, i.e., $C^1(L,\KK)$ does not contain a non--trivial
$L_0$--invariant element. Consequently, Prop.~2.1 shows that $H^2(L,\KK)$ can
be identified with the space of all $L_0$--invariant 2--cochains. Thus we are
going to determine this space.

For $n \ge 3$\,, this is easy. Let us first recall that every $L_0$--invariant
bilinear form on $L_0$ is symmetric. Using this fact, a short look at the
representations of $L_0$ carried by $L_0$ and $L_{\pm 1}$ shows that there
exists, up to the normalization, a unique non--zero super--skew--symmetric
$L_0$--invariant bilinear form on $L$\,. As we already know, $g$ has these
properties, moreover, $g$ is a 2--cocycle. Thus we conclude, for $n \ge 3$\,,
that $H^2(L,\KK)$ is one--dimensional and that \reff{slcov} is a universal
covering of $L = sl(n|n)/\KK \cd I_{2n}$\,.

Due to the peculiarities of the algebra $sl(2|2)/\KK \cd I_4$\,, the case
$n = 2$ is more interesting. In principle, we could proceed in this case as
above. However, we prefer to use a different realization of
$sl(2|2)/\KK \cd I_4$ which makes the special properties of this algebra
manifest (and which has already been used in the description of the super
derivations of this algebra \cite{S}).

Let $V$ be a two--dimensional vector space, let $\psi$ be a non--degenerate
skew--symmetric bilinear form on $V$, and let $sp(\psi)$ (equal to $sl(V)$)
be the Lie algebra of all linear maps of $V$ into itself leaving the form
$\psi$ invariant. Define the bilinear map
 \be  P : V \ti V \lra sp(\psi)  \ee
by
 \be  P(x,y)\,z = \psi(x,z)\,y + \psi(y,z)\,x  \ee
for all $x,y,z \in V$. We set
 \be  \Lzb = sp(\psi) \ti sp(\psi) \;\;\;,\;\;\; \Lob = V \ot V \ot V  \ee
and regard $\Lob$ as an $\Lzb$--module, where the first (resp. second)
factor of $\Lzb$ acts in the natural way on the first (resp. second)
tensorial factor of $\Lob$\,. For any $Q \in sp(\psi)$, let $Q_1$ and $Q_2$
be the elements of $\Lzb$ defined by
 \be  Q_1 = (Q,0) \;\;\;,\;\;\; Q_2 = (0,Q) \,.  \ee
Using these data, the $\ZZ_2$--graded vector space
 \be  L = \Lzb \oplus \Lob  \ee
becomes a Lie superalgebra, if the bracket of two elements of $\Lob$ is
chosen such that
 \be  \la x \ot y \ot z,x' \ot y' \ot z' \ra
 = c\,P_1(x,x')\psi(y,y')\psi(z,z') \,-\, c\,\psi(x,x')P_2(y,y')\psi(z,z')  \ee
for all $x,y,\ldots,z' \in V$, where $c \in \KK$ is an arbitrary non--zero
constant.

This Lie superalgebra is isomorphic to $sl(2|2)/\KK \cd I_4$\,. More precisely,
one can easily show the following: \\
Let $(v_{+},v_{-})$ be any basis of $V$. If we define
 \be  L_0 = \Lzb \;\;\;,\;\;\; L_{\pm 1} = V \ot V \ot \KK\,v_{\pm} \,, \ee
then
 \be  L = L_{-1} \oplus L_0 \oplus L_1  \ee
is a consistent $\ZZ$--gradation of $L$\,, and $L$ and $sl(2|2)/\KK \cd I_4$
are isomorphic as $\ZZ$--graded Lie superalgebras.

The arbitrariness in the choice of $v_{\pm}$ reflects the peculiarity of
$sl(2|2)/\KK \cd I_4$ we have been alluding to: This algebra admits a large
group of outer automorphisms. In fact, let $s$ be an arbitrary element of
the group $Sp(\psi)$ (equal to $SL(V)$) and let $\tilde{s}$ be the linear
map of $L$ into itself, which induces the identity map on $\Lzb$\,, and on
$\Lob$ is given by
 \be  \tilde{s}(x \ot y \ot z) = x \ot y \ot s(z)  \ee
for all $x,y,z \in V$. Then $\tilde{s}$ is an automorphism of the Lie
superalgebra $L$\,.

It is easy to determine the $\Lzb$--invariant super--skew--symmetric bilinear
forms on $L$\,. Let $g$ be such a form. For the same reasons as above,
$g(A,B)$ has to vanish if at least one of the two homogeneous elements
$A,B \in L$ is even. On the other hand, the restriction of $g$ to
$\Lob \ti \Lob$ is given by
 \be g(x \ot y \ot z,x' \ot y' \ot z') = \psi(x,x')\psi(y,y')\vp(z,z') \ee
for all $x,y,\ldots,z' \in V$, where $\vp$ is an arbitrary symmetric bilinear
form on $V$. Conversely, it is not difficult to check that every bilinear form
$g$ with these properties is, in fact, an $\Lzb$--invariant 2--cocycle.

Obviously, the map assigning to any symmetric bilinear form $\vp$ on $V$ the
$\Lzb$--invariant 2--cocycle specified above is a vector space isomorphism of
the spaces under consideration. Consequently, we have
 \be \mbox{dim}\,H^2(sl(2|2)/\KK \cd I_4\,,\KK) = 3 \,.  \ee

The universal covering of $L$ can now be constructed along the lines specified
around Eq.~\reff{univcocy}. Actually, this is not really necessary, since this
algebra is known: It is obtained if one sets out to construct the Serre
presentation of $sl(2|2)$ and forgets about the so--called supplementary
relations. More precisely, the universal covering of $sl(2|2)/\KK \cd I_4$
is the Lie superalgebra $s(2,2)$ considered in Ref.~\cite{S.pq}. In fact,
there is a natural homomorphism of the Lie superalgebra $s(2,2)$ onto the
Lie superalgebra $sl(2|2)/\KK \cd I_4$\, whose kernel is equal to the
three--dimensional center of $s(2,2)$, moreover, it is known that $s(2,2)$
is perfect. Since the universal covering of $sl(2|2)/\KK \cd I_4$ must cover
$s(2,2)$ as well (see Remark 4.1), it follows that this universal covering
must be isomorphic to $s(2,2)$. Actually, it is not difficult to construct
this isomorphism explicitly. According to Remark 4.1, $s(2,2)$ is also a
universal covering of $sl(2|2)$. (Note that this is at variance with
Thm.~2.6.1 of Ref.~\cite{Fu}, but it is consistent with Ref.~\cite{FL}. On the
other hand, in the latter reference, the definition of the two basic
2--cocycles of $sl(2|2)$ is not correct, and it is easy to see that, contrary
to what is stated, we have $H^3(L,\KK) \neq \Ze$ for both $L = sl(2|2)$ and
$L = sl(2|2)/\KK \cd I_4$\,. We think that these discrepancies and mistakes
are due to some trivial slips or misprints.)
\vspace{5ex}

%
\noindent
{\Large \bf
V. Cohomology of the Lie superalgebra $sl(1|2)$ \\[-0.5ex]}
\setcounter{section}{5} \setcounter{equation}{0}

\noindent
Throughout the present section we shall assume that the base field $\KK$ is
{\em algebraically closed.} Let $U(sl(1|2))$ denote the enveloping algebra
of $sl(1|2)$. Our main goal is to show that
 \be  H^2(sl(1|2),U(sl(1|2))) = \Ze \,. \label{chslvan} \ee
In order to prove this, we shall first determine $H^1(sl(1|2),V)$ and
$H^2(sl(1|2),V)$, for all finite--dimensional simple graded $sl(1|2)$--modules
$V$. Actually, we only need to know those modules $V$ for which these
cohomology groups are non--trivial. As explained in the introduction, this
information could be extracted from Ref.~\cite{Ta}, but we think it is
worthwhile to present our own calculations.

As an abbreviation, we write $L$ instead of $sl(1|2)$. We are going to use
the notation and results of Ref.~\cite{SNR}. Let us first recall some of the
most relevant facts. The basis of $L$ we are going to use consists of four
even elements $\Qpm$\,, $Q_3$\,, $B$ and four odd elements $\Vpm$ and $\Wpm$\,.
The elements $\Qpm$\,, $Q_3$ span an $sl(2)$ subalgebra of $\Lzb$\,, they are
normalized such that
 \be  \la Q_{+},Q_{-}\ra = 2 Q_3 \,.  \ee
The one--dimensional subspace $\KK \cd B$ is the center of $\Lzb$\,, and the
two--dimensional subspace
 \be  H = \KK \cd B + \KK \cd Q_3  \ee
is a Cartan subalgebra of $L$\,. Thus any linear form $\lam$ on $H$ (in
particular, any weight of an $L$--module) can and will be identified with the
pair $(b,q)$ of values which $\lam$ takes on $B$ and $Q_3$\,, respectively.
For example, $\Vpm$ and $\Wpm$ are root vectors corresponding to the roots
$(\hal,\pm\hal)$ and $(-\hal,\pm\hal)$, respectively. They are normalized such
that $(\Vpm)$ and $(\Wpm)$ are standard $sl(2)$--doublets and satisfy
 \bea      \la \Vpm ,\Wpm\ra \eqq \pm\Qpm  \\
        \la \Vpm ,W_{\mp}\ra \eqq -Q_3 \pm B \,.  \eea
The reader who prefers a matrix realization may set
 \be Q_{+} = E_{23} \;\;,\;\; Q_{-} = E_{32} \;\;,\;\;
     Q_3 = \ohal(E_{22} - E_{33}) \;\;,\;\;
       B = - \ohal(2E_{11} + E_{22} + E_{33})  \ee
 \be V_{+} = E_{21} \;\;,\;\; V_{-} = E_{31} \;\;,\;\;
     W_{+} = E_{13} \;\;,\;\; W_{-} = - E_{12} \,,  \ee
where the $E_{ij}$ are the usual basic $3 \ti 3$ matrices. (The elements are
chosen such that the $\ZZ$--gradation of $sl(1|2)$ as specified in
Ref.~\cite{S} and the $\ZZ$--gradation given by Eq.~\reff{zgrad} coincide.)

For later use we note that the linear map
 \be  \om : L \lra L \label{ombeg} \ee
defined by
 \be  \om(\Qpm) = \Qpm \;\;\;,\;\;\; \om(Q_3) = Q_3 \;\;\;,\;\;\;
                                                           \om(B) = -B  \ee
 \be  \om(\Vpm) = \Wpm \;\;\;,\;\;\; \om(\Wpm) = \Vpm \label{omend} \ee
is an automorphism of the Lie superalgebra $L$\,.

The simple roots of $L$ with respect to $H$ are chosen such that $V_{+}$ and
$W_{+}$ are the associated root vectors. Thus both of the simple roots are
odd, i.e., we do not use a so--called distinguished basis of the root system.
For $sl(1|2)$, this choice has certain advantages, for example, it is
compatible with the automorphism $\om$\,.

For any weight $(b,q)$\,; $b,q \in \KK$\,, and any $\ga \in \ZZ_2$\,, there
exists (up to isomorphism) a unique simple graded highest weight $L$--module
$V(b,q,\ga)$ such that the highest weight vectors (which are proportional to
each other) have weight $(b,q)$ and $\ZZ_2$--degree $\ga$\,. Of course,
$V(b,q,\ga)$ and $V(b,q,\ga + \ob)$ can be transformed into each other
by a shift of the $\ZZ_2$--gradation (see Eq.~\reff{grshift}). The module
$V(b,q,\ga)$ is finite--dimensional if and only if $b = q = 0$ or
$q \in \hal\NS$\,, and any finite--dimensional simple graded $L$--module is
of this type. Moreover, it is atypical \cite{Ka.bo} if and only if
$b = \pm q$\,, and its dimension is equal to $8q$ (resp.~$4q+1$) in the
typical (resp.~atypical) case.

The Casimir elements of $L$ are known, and the eigenvalues of the Casimir
operators in a graded highest weight module have been calculated in terms of
a generating function \cite{S.ce,S.cec,S.ev} (see also Example 4 of Sec.~III).
In particular, the image of the generalized Harish--Chandra isomorphism is
known. It can be regarded as an algebra of polynomials in $b,q$ and, as such,
consists of all polynomials of the form $a + (b^2 - q^2)\,P$\,, where
$a \in \KK$ is a constant and $P$ is a polynomial in $b$ and $(b^2 - q^2)$.
This implies that, for a graded highest weight module with highest weight
$(b,q)$, all Casimir operators without a constant term are equal to zero, if
and only if $b = \pm q$\,.

Next we recall that $L$ has a consistent $\ZZ$--gradation, which is fixed by
the requirement that the elements $\Qpm$\,, $Q_3$\,, and $B$ are homogeneous
of degree $0$\,, the elements $\Vpm$ are homogeneous of degree $1$\,, and the
elements $\Wpm$ are homogeneous of degree $-1$\,. Thus we can apply our
general theory both in the $\ZZ_2$--graded and in the $\ZZ$--graded sense.
Actually, this essentially amounts to a matter of formulation (at least in
the present article). This is because, for all $L$--modules $V$ for which a
$\ZZ$--gradation will be of interest, this gradation is simply given by the
action of $B$\,, as follows:
 \be  V_r = \{x \in V \,|\, 2B \cd x = r\,x \} \label{zgrad} \ee
for all $r \in \ZZ$\,. In particular, this is true for $L$ itself and for the
enveloping algebra $U(L)$ (both endowed with the adjoint action of $L$\,).

As a final preparatory remark we notice that there is a natural embedding of
the Lie superalgebra $osp(1|2)$ into $sl(1|2)$. In fact, if we define
 \be  \Upm = \ohal(\Vpm + \Wpm) \,,  \ee
the elements $\Qpm$\,, $Q_3$\,, $\Upm$ satisfy the commutation relations of
$osp(1|2)$ as specified in Ref.~\cite{SNR}. In particular, we have
 \bea     \la \Upm,\Upm \ra \eqq \pm\ohal\Qpm  \\
       \la \Upm,U_{\mp} \ra \eqq -\ohal Q_3 \,.  \eea
The $\ZZ_2$--graded subalgebra of $L$ spanned by these elements will be denoted
by $G$\,, by construction, it is isomorphic to $osp(1|2)$.

We recall that the finite--dimensional graded $G$--modules are semi--simple.
The finite--dimensional simple graded $G$--modules can be labelled by a pair
$(q,\ga)$, with $q \in \hal\NN$ and $\ga \in \ZZ_2$\,: By definition, a
highest weight vector $v$ of the corresponding module is homogeneous of
degree $\ga$ and satisfies $Q_3 \cd v = q\,v$\,. We denote this module (which
is fixed up to isomorphism) by $U(q,\ga)$, its dimension is equal to $4q+1$\,.

Of course, every graded $L$--module can also be regarded as a graded
$G$--module. For example, choose $b \in \KK$\,, $q \in \hal\NS$\,,
$\ga \in \ZZ_2$\,, and regard $V(b,q,\ga)$ as a $G$--module. If
$b \neq \pm q$\,, this module decomposes into the direct sum of two simple
graded $G$--modules, isomorphic to $U(q,\ga)$ and $U(q - \hal,\ga + \ob)$,
respectively. On the other hand, the $G$--module $V(\pm q,q,\ga)$ is simple
and isomorphic to $U(q,\ga)$. (Obviously, the latter statement remains true
for $q = 0$\,.)

In particular, consider the adjoint $L$--module $L$\,. Regarded as a
$G$--module, it decomposes into two $G$--submodules: One of them is $G$
itself, the other is spanned by $\Xpm$ and $B$\,, with
 \be  \Xpm = \ohal(\Vpm - \Wpm) \,.  \ee
The elements $\Xpm$ have been chosen such that $\Xpm$ and $B$ form a basis of
$U(\hal,\ob)$ as specified in Ref.~\cite{SNR}. We remark that $\Qpm$\,,
$Q_3$\,, $i \Xpm$ also satisfy the commutation relations of $osp(1|2)$\,.

After these preliminaries we are ready to tackle our main problem: We want to
determine the cohomology groups $H^1(L,V)$ and $H^2(L,V)$, with $V$ a
finite--dimensional simple graded $L$--module. In view of the preceding
remarks, we deduce from Prop.~2.2 that $H^n(L,V) = \Ze$ for all $n$ if $V$ is
typical. Consequently, we only have to consider the cases where
$V = V(\pm q,q,\ga)$, with $q \in \hal\NN$ and $\ga \in \ZZ_2$\,. Actually,
if $\vr$ is the representation afforded by $V(q,q,\ga)$, then $\vr \comp \om$
is (equivalent to) the representation afforded by $V(-q,q,\ga)$ (where the
automorphism $\om$ has been defined by the Eqs.~\reff{ombeg}\,--\,\reff{omend}).
Recalling the isomorphism \reff{omcohom} we may assume, therefore, that
$V = V(q,q,\ga)$. On the other hand, the bijection \reff{hshift} shows that
it is sufficient to consider only one of the two choices for $\ga$\,.
According to Eq.~\reff{zgrad}, we shall choose $\ga = \tqb$ (where
$\bar{r} = r + 2\ZZ \in \ZZ_2$ denotes the residue class of the integer $r$
modulo 2\,), and assume that $V(q,q,\tqb)$ is endowed with the
consistent $\ZZ$--gradation defined by Eq.~\reff{zgrad}. To simplify the
notation, we write $V(q)$ instead of $V(q,q,\tqb)$. Then we have to
determine the cohomology groups $H^n(L,V(q))$, with $n \in \{1,2\}$.

The case $n = 1$ is simple. We present it in some detail to explain our
approach. According to Prop.~2.1, any 1--cocycle of $L$ with values in $V(q)$
is cohomologous to a $G$--invariant 1--cocycle. Regarded as graded
$G$--modules, $L$ is isomorphic to $U(1,\zb) \op U(\hal,\ob)$, and $V(q)$ is
isomorphic to $U(q,\tqb)$. But a non--zero $G$--invariant linear map of
a $G$--module $U(p,\tpb)$ into a $G$--module $U(q,\tqb)$ (with
$p,q \in \hal\NN$\,) exists if and only if $p = q$\,, and such a map must be
even. Consequently, a non--zero $G$--invariant linear map $g : L \rar V(q)$
exists if and only if $q \in \{\hal,1\}$. Moreover, if this condition is
satisfied, a map of the type in question is even and uniquely determined up
to a scalar factor.

On the other hand, let $V$ be any graded $L$--module, and let $g : L \rar V$
be a $G$--invariant 1--cocycle which is homogeneous of degree $\zb$\,. Then
the $G$--invariance and the cocycle condition imply that, for all homogeneous
elements $X \in G$ and $Y \in L$\,,
 \be  X \cd g(Y) = g(\la X,Y\ra) = X \cd g(Y) - \ve(\xi,\eta)\,Y \cd g(X) \ee
and hence that
 \be  Y \cd g(X) = 0 \,.  \ee
In particular, we conclude that $g(\la X,X'\ra) = 0$ for all $X,X' \in G$\,,
which implies that
 \be  g(X) = 0  \ee
for all $X \in G$\,.

Thus all we have to do is to determine $H^1(L,V(\hal))$. To construct a
non--zero $G$--module homomorphism of $L$ into $V(\hal)$ we introduce a
suitable basis in $V(\hal)$. According to Ref.~\cite{SNR} there exists a
basis $(e_{+},e_0,e_{-})$ of $V(\hal)$ such that $\epm$ are odd and $e_0$
is even, such that $\epm$ and $e_0$ are weight vectors corresponding to the
weights $(\hal,\pm\hal)$ and $(1,0)$, respectively, and such that
 \be  \Vpm \cd \epm = 0 \;\;\;,\;\;\; \Vpm \cd e_{\mp} = \mp e_0  \ee
 \be          \Wpm \cd \epm = \Wpm \cd e_{\mp} = 0  \ee
 \be  \;\Vpm \cd e_0 = 0 \;\;\;,\;\;\; \Wpm \cd e_0 = -\epm \,. \ee
(Obviously, this basis is uniquely fixed up to a non--zero overall factor.)
The action of $\Qpm$ on the basis vectors can be derived from the commutation
relations and the formulae above, it is such that $(\epm)$ is a standard
$sl(2)$--doublet and $e_0$ an $sl(2)$--singlet. It follows that $\epm$ and
$e_0$ also form a basis of the $G$--module $V(\hal) \simeq U(\hal,\ob)$ as
described in Ref.~\cite{SNR}. Consequently, the linear map
 \be  g : L \lra V(\ohal) \label{gbeg} \ee
defined by
 \be  g(X) = 0 \;\;\;\mbox{if}\;\;\; X \in G  \ee
 \be  g(B) = e_0  \ee
 \be  g(\Vpm - \Wpm) = 2\epm \label{gend} \ee
is a homomorphism of graded $G$--modules. Since
 \be  g(\Vpm + \Wpm) = g(2\Upm) = 0 \,,  \ee
it follows that
 \be  g(\Vpm) = -g(\Wpm) = \epm \,.  \ee
It is now easy to check that $g$ is a 1--cocycle. Moreover, since $V(\hal)$
does not contain a non--zero $G$--invariant element, $g$ is not a coboundary.
Summarizing, we have shown that
 \be  H^1(L,V(q)) = \Ze \;\;\;\mbox{if}\;\;\; q \neq \ohal  \ee
 \be  \mbox{dim}\,H^1(L,V(\ohal)) = 1 \,.  \ee
The fact that $V(\hal)$ is (essentially) the sole finite--dimensional simple
graded $L$--module $V$ for which $H^1(L,V)$ is non--trivial has already been
mentioned in the {\em Additional remarks} at the end of Ref.~\cite{Ka}.

The cocycle $g$ can be replaced by a simpler one. Let us regard $L$ as a
consistently $\ZZ$--graded Lie superalgebra and $V(\hal)$ as a $\ZZ$--graded
$L$--module (recall that the $\ZZ$--gradations are chosen according to
Eq.~\reff{zgrad}). Then $g$ can be decomposed into its $\ZZ$--homogeneous
components,
 \be  g = g_0 + g_2 \,, \ee
where $g_r$\,, $r \in \{0,2\}$, is homogeneous of degree $r$\,. Explicitly,
the non--vanishing values of $g_0$ and $g_2$ on the basis vectors $\Qpm$\,,
$Q_3$\,, $B$\,, $\Vpm$\,, $\Wpm$ are given by
 \be  g_{0}(\Vpm) = \epm  \ee
 \be  g_{2}(B) = e_0 \;\;\;,\;\;\; g_{2}(\Wpm) = -\epm \,.  \ee
Both $g_0$ and $g_2$ are 1--cocycles. Hence exactly one of them is not a
coboundary. Obviously, $g_0$ cannot be a coboundary, for $V(\hal)$ does not
contain a non--zero $\ZZ$--homogeneous element of degree 0\,. Consequently,
$g_2$ must be a coboundary. In fact, we have
 \be  g_{2}(X) = X \cd e_0 \label{gtwo} \ee
for all $X \in L$\,.

The 1--cocycle $g_0$ can be used to construct higher order cocycles, as
follows. Choose any integer $n \ge 1$\,. Using the product of cocycles
defined in Eq.~\reff{prodall} and the associativity of this product (see
Eq.~\reff{prodass}), we can construct the $n$--cocycle
$g_{0}^{\od n} = g_0 \od\ldots\od g_0$ ($n$ factors $g_0$\,) with values in
$V(\hal)^{\ot n}$, which is given by
 \be g_{0}^{\od n}(X_1,\ldots,X_n) = \sum_{\pi \in \Sn}
                     \!\mbox{sgn}(\pi)\,\ve_n(\pi;\xi_1,\ldots,\xi_n)
                     \,g_0(X_{\pi(1)}) \ot\ldots\ot g_0(X_{\pi(n)})  \ee
for all homogeneous elements $X_i \in L$\,. Obviously, $g_{0}^{\od n}$ is
non--zero and takes its values in the graded $L$--submodule $W(n)$ of
super--skew--symmetric tensors in $V(\hal)^{\ot n}$.

Let $f_n : L^n \rar W(n)$ denote the $n$--cocycle defined by $g_{0}^{\od n}$.
Of course, $f_n$ is $\ZZ$--homogeneous of degree zero. Since the degree
of a non--zero $\ZZ$--homogeneous element of $W(n)$ is at least equal to
$n$\,, there does not exist a non--zero $(n-1)$--linear map of $L^{n-1}$
into $W(n)$, which is $\ZZ$--homogeneous of degree zero. Consequently, $f_n$
is not a coboundary. On the other hand, it is easy to see that $W(n)$, even
when regarded as a $\ZZ$--graded $L$--module, is isomorphic to
$V(\frac{n}{2})$. Thus we have proved that
 \be  H^n(L,V(\nhal)) \neq \Ze \;\;\;\mbox{for all}\;\;\; n \ge 1 \,.
                                                         \label{hnneqz} \ee

\noindent
{\it Remark 5.1.} Our results show that Prop.~2.2 exactly describes the
finite--dimensional simple graded $L$--modules with a non--trivial cohomology,
and the same holds for the $osp(1|2n)$ algebras. It would be interesting to
know whether this is also true for other simple Lie superalgebras.

Let us next determine the cohomology groups $H^2(L,V(q))$. Since we are going
to proceed as for $H^1(L,V(q))$, it should suffice to give the main
intermediate steps. For any graded $L$--module $V$, the graded $L$--module
$C^2(L,V)$ is canonically isomorphic to $\Lg(L \wepf L,V)$ (see App.~A). Thus
we have to determine the structure of $L \wepf L$\,. It is known that
this module is isomorphic to the graded submodule of super--skew--symmetric
tensors in $L \ot L$\,. This submodule is semi--simple (whereas $L \ot L$ is
not), and its simple submodules can be determined, with the result that
 \be  L \wepf L \,\simeq\, V(\ohal,\thhal,\ob) \op V(-\ohal,\thhal,\ob)
                                                     \op V(0,1,\zb) \,.  \ee
This implies that, regarded as a graded $G$--module, $L \wepf L$ is isomorphic
to
 \be  2\,U(\thhal,\ob) \op 3\,U(1,\zb) \op U(\ohal,\ob) \,.  \ee
Recalling that the graded $G$--module $V(q)$ is isomorphic to $U(q,\tqb)$,
we conclude as before that $H^2(L,V(q))$ is trivial unless
$q \in \{\hal,1,\thhal\}$.

Let us consider the case $q = \hal$\,. Visibly, there exists, up to a scalar
factor, a unique non--zero $G$--invariant linear map $\bar{f}$ of $L \wepf L$
into $V(\hal)$, and this map is even and vanishes on the submodule $W$
isomorphic to
 \be  W \,\simeq\, V(\ohal,\thhal,\ob) \op V(-\ohal,\thhal,\ob) \,.  \ee
Let $f \in C^2(L,V(\hal))$ be the associated bilinear map. We are going to
show that $f$ is not a cocycle. According to Prop.~2.1, this will imply that
 \be  H^2(L,V(\ohal)) = \Ze \,.  \ee
In fact, the product map $\la\;,\;\ra$ of $L$ induces an $L$--module
homomorphism of $L \wepf L$ onto $L$\,, whose kernel is equal to $W$.
Consequently, there exists a unique homomorphism of graded $G$--modules
 \be  g' : L \lra V(\ohal)  \ee
such that
 \be  f(X,Y) = g'(\la X,Y\ra)  \ee
for all $X,Y \in L$\,. But we know that there exists, up to a scalar factor,
a unique homomorphism of graded $G$--modules of $L$ into $V(\hal)$, namely,
the map $g$ defined by the Eqs.~\reff{gbeg}\,--\,\reff{gend}. Thus, by
multiplying $f$ by a suitable scalar factor, we may assume that
 \be  f(X,Y) = g(\la X,Y\ra)  \ee
for all $X,Y \in L$\,. If $f$ would be a 2--cocycle, its $\ZZ$--homogeneous
components $f_r$ given by
 \be  f_r(X,Y) = g_r(\la X,Y\ra) \;\;\;,\;\;\; r \in \{0,2\}  \ee
would be 2--cocycles as well. However, using Eq.~\reff{gtwo}, it is easy to
see that $f_2$ is not a cocycle.

To proceed, we remark that the condition $q \in \{\hal,1,\thhal\}$ derived
above can be improved by regarding $L \wepf L$ and $V(q)$ as $\Lzb$--modules.
The same type of argument as before then shows that $H^2(L,V(q))$ is trivial
unless $q \in \{0,\hal,1\}$, moreover, there exists, up to a scalar factor, a
unique $\Lzb$--invariant linear map of $L \wepf L$ into $V(1)$. Since (again
by Prop.~2.1) any 2--cocycle on $L$ with values in $V(1)$ is cohomologous to
an $\Lzb$--invariant cocycle, we conclude from Eq.~\reff{hnneqz} that
$H^2(L,V(1))$ is one--dimensional. Summarizing, we have shown that
 \be  H^2(L,V(q)) = \Ze \;\;\;\mbox{if}\;\;\; q \neq 1 \label{chvqvan} \ee
 \be  \mbox{dim}\,H^2(L,V(1)) = 1 \,.  \ee

As an application of the preceding result, let us now prove that $H^2(L,U(L))$
is trivial. The following lemma contains the pertinent information on the
$L$--module $U(L)$.

\noindent
{\bf Lemma 5.1} \\
{\it Any finite--dimensional atypical simple graded subquotient of the adjoint
$L$--module $U(L)$ is isomorphic to one of the three graded $L$--modules
$V(0,0,\zb)$ and $V(\pm\hal,\hal,\ob)$. }

\noindent
Proof \\
Consider two graded submodules $U' \subset U$ of $U(L)$ such that $U/U'$ is a
finite--dimensional atypical simple graded $L$--module. Then $U/U'$ is
isomorphic to one of the modules $V(\pm q,q,\ga)$, with $q \in \NN$ and
$\ga \in \ZZ_2$\,, and we have to show that $q \in \{0,\hal\}$ and
$\ga = \tqb$\,.

Let $\ad$ denote the adjoint representation of $L$ in $U(L)$. Then $\ad B$ is
diagonalizable, and the sole eigenvalues of this operator are $0$\,,
$\pm\hal$\,, and $\pm 1$\,. Since, for $q \geq \hal$\,, the module
$V(\pm q,q,\ga)$ contains non--zero weight vectors with weight
$(\pm q \pm \hal,q - \hal)$, the cases $q \geq 1$ are not possible. On the
other hand, an eigenvector of $\ad B$ is odd if and only if the corresponding
eigenvalue is equal to $\pm\hal$\,, and it is even otherwise. This implies
that $\ga = \tqb$ and proves the lemma.

\noindent
{\it Remark 5.2.} As shown in App.~B, the adjoint $L$--module $U(L)$ really
has graded subquotients isomorphic to $V(0,0,\zb)$ and $V(\pm\hal,\hal,\ob)$.

\noindent
{\it Remark 5.3.} A decomposition of the adjoint $L$--module $U(L)$ into
indecomposable submodules has been constructed in a recent paper by Benamor
\cite{Be}. (We are grateful to the referee for drawing our attention to this
article.) This decomposition is rather complicated, thus we think it would
not be adequate to invoke that reference to prove the simple lemma above.

It is now easy to see that $H^2(L,U(L)) = \Ze$. The image of any 2--cochain
of $L$ with values in $U(L)$ is contained in a finite--dimensional graded
$L$--submodule of $U(L)$ (use the canonical filtration of $U(L)$). Thus it
is sufficient to show that $H^2(L,V)$ is trivial, for any finite--dimensional
graded $L$--submodule $V$ of $U(L)$. To prove this we choose an increasing
sequence of graded $L$--submodules of $U(L)$,
 \be  \Ze = V_0 \subset V_1 \subset V_2 \subset \ldots \subset V_n = V \,,
                                                              \label{jhs} \ee
such that $V_r/V_{r-1}$ is a simple graded $L$--module, for $1 \le r \le n$\,.
(Up to the opposite ordering, \reff{jhs} is a Jordan--H\"older sequence of the
graded $L$--module $V$.) We show by induction on $r$ that
 \be  H^2(L,V_r) = \Ze \,. \label{hrvan} \ee
The case $r = 0$ is trivial. Now suppose that Eq.~\reff{hrvan} is true for
some $r \in \{0,1,\ldots,n-1\}$. We have a short exact sequence of graded
$L$--modules:
 \be \Ze \lra V_r \lra V_{r+1} \lra V_{r+1}/V_r \lra \Ze \,.\label{shexseq}\ee
Consider the following part of the exact cohomology sequence \reff{cohomseq}
associated to \reff{shexseq}:
 \be  \ldots\,\lra H^2(L,V_r) \lra H^2(L,V_{r+1}) \lra H^2(L,V_{r+1}/V_r)
                                                       \lra\,\ldots  \;.  \ee
By assumption, $H^2(L,V_r)$ is trivial, and Lemma 5.1, combined with
Eq.~\reff{chvqvan}, shows that $H^2(L,V_{r+1}/V_r)$ is trivial as well. Thus,
$H^2(L,V_{r+1})$ must be trivial.

Summarizing, we have shown that
 \be  H^2(L,U(L)) = \Ze \,, \label{van} \ee
i.e., we have proved Eq.~\reff{chslvan}. On the other hand, we shall see
in App.~B that
 \be  H^1(L,U(L)) \neq \Ze \,,  \ee
in contrast to what is known for semi--simple Lie algebras.
\vspace{5ex}

%
\noindent
{\Large \bf
VI. Discussion \\[-0.5ex]} \setcounter{section}{6} \setcounter{equation}{0}

\noindent
In the present work we have taken some exploratory steps towards a better
understanding of the cohomology of Lie superalgebras and their generalizations.
It is hardly surprising that the basic definitions and formal techniques known
from the cohomology of Lie algebras can immediately be generalized to the
general graded setting. In particular, this applies to Garland's theory of
universal central extensions of Lie algebras. The picture changes if we
try to actually calculate the cohomology groups of Lie superalgebras,
especially, when the coefficients are non--trivial. Our main tools in this
project were two simple propositions proved in Sec.~II.

As for many questions from the theory of Lie superalgebras, the $osp(1|2n)$
algebras behave very much like simple Lie algebras. On the other hand, already
for $sl(1|2)$ the situation changes drastically. For any integer $n \ge 0$\,,
there exists a finite--dimensional simple graded $sl(1|2)$--module $V$ such
that $H^n(sl(1|2),V)$ is non--trivial. Moreover, in App.~B we are going to
show that
 \be  H^1(sl(1|2),U(sl(1|2))) \neq \Ze \,.  \ee
Thus it is remarkable that, nevertheless,
 \be  H^2(sl(1|2),U(sl(1|2))) = \Ze \,.  \ee
Actually, one of the main reasons for our present investigation was to try
to prove the analogous result for as many simple Lie superalgebras $L$ as
possible, from which it follows that the associative superalgebra $U(L)$ does
not admit of non--trivial formal deformations in the sense of Gerstenhaber
\cite{Ge}. To see this, all we have to do is to transcribe the corresponding
discussion of Ref.~\cite{Kas} to the present setting. Then Eq.~\reff{van}
implies that the Hochschild 2--cohomology (in the super sense) of the
associative superalgebra $U(L)$ with values in $U(L)$ (considered as a graded
$U(L)$--bimodule) vanishes. In this connection, we recall that there is
an isomorphism relating the cohomology of a Lie algebra $L$ with values in
$U(L)$ to the Hochschild cohomology of the associative algebra $U(L)$ (see
Theorem 5.1 of chapter XIII in Ref.~\cite{CaE}). It would be interesting to
know whether an analogous isomorphism exists in the general graded setting.

The methods used in our investigation are completely elementary, and already
for the 3--cohomology of $sl(1|2)$ the calculations become quite extensive.
On the other hand, in the meantime we have shown that $H^2(L,U(L)) = \Ze$
for $L = sl(m|1)$, $m \geq 3$\,, but up to now we were not able to prove or
disprove this for $L = sl(3|2)$. Thus it might turn out that some more profound
techniques must be used both in the study of higher order cohomology groups
and in the case of more complicated Lie superalgebras.

\vspace{2ex}

\noindent
{\bf Acknowledgement} \\
We should like to thank C. Gruson for sending us a copy of the galley proofs
of her article \cite{Gr} and for drawing our attention to Ref.~\cite{Ta}.
The present work was initiated during a visit of the first--named author to
the Department of Pure Mathematics of the University of Adelaide. The kind
invitation by the second--named author and the hospitality extended to the
first--named author, both in the Mathematics and in the Physics Department,
are gratefully acknowledged.

\vspace{8ex}  

\noindent
{\LARGE\bf Appendix}\\[-6ex]

\begin{appendix}

\sect{The $\ve$--exterior algebra of a graded $L$--module}
We use the notation introduced at the beginning of Sec.~II. Let $V$ be a graded
$L$--module and let
 \be  T(V) = \bop_{n \in \ZZ}T_n(V)  \ee
be the tensor algebra of $V$ \cite{S.ma}. It is known that $T(V)$ is a
$\ZZ \ti \Ga$--graded associative algebra with a unit element and that
 \bea  T_n(V) \eqq \Ze \;\;\;\mbox{if}\;\;\; n \le -1 \\
       T_0(V) \eqq \KK \\
       T_n(V) \eqq
 \underbrace{V \ot \ldots \ot V}_{n\;\mbox{{\scriptsize factors}}}
                                \;\;\;\mbox{if}\;\;\; n \ge 1 \,. \eea
Moreover, $T(V)$ is a graded $L$--module, and $L$ acts on $T(V)$ by
$\ve$--derivations. Note also that each $T_n(V)$ is a graded $L$--submodule
of $T(V)$.

Now let $J(V,\ve)$ denote the two--sided ideal of $T(V)$, generated by the
elements of the form
 \be  x \ot y + \ve(\xi,\eta)\,y \ot x \,,  \ee
where $x$ and $y$ are homogeneous elements of $V$. Then the quotient algebra
 \be  \Wepf V = T(V)/J(V,\ve)  \ee
is called the {\it $\ve$--exterior (or $\ve$--Grassmann) algebra} of $V$. The
multiplication in $\Wept V$ will be denoted by $\wepe$\,.

Since $J(V,\ve)$ is a $\ZZ \ti \Ga$--graded ideal, the algebra $\Wept V$
inherits from $T(V)$ a canonical $\ZZ \ti \Ga$--gradation. In particular, we
have
 \be  \Wepf V = \bop_{n \in \ZZ}\Wepfn V \,,  \ee
where $\Weptn V$ is the canonical image of $T_n(V)$ in $\Wept V$. It follows
that
 \bea  \Wepfn V \eqq \Ze \;\;\;\mbox{if}\;\;\; n \le -1 \\
       \Wepfz V \eqq \KK  \\
       \Wepfo V \eqq V \,.  \eea
For $n \ge 2$ we also use the notation
 \be  \Wepfn V = V \wepf \ldots \wepf V \;\;\;(n\;\mbox{factors}) \,.  \ee
On the other hand, $J(V,\ve)$ is a $\Ga$--graded $L$--submodule of $T(V)$, and
hence $\Wept V$ inherits from $T(V)$ the structure of a $\Ga$--graded
$L$--module. Then $L$ acts on $\Wept V$ by $\ve$--derivations, and each
$\Weptn V$ is a $\Ga$--graded $L$--submodule of $\Wept V$. 

The theory of $\ve$--exterior algebras can easily be reduced to the
theory of $\ve'$--symmetric algebras \cite{S.ma}, with $\ve'$ a commutation
factor on a suitable Abelian group $\Ga'$. In fact, set $\Ga' = \ZZ \ti \Ga$
and define the commutation factor $\ve'$ on $\Ga'$ by
 \be  \ve'((r,\al),(s,\bet)) = (-1)^{rs}\ve(\al,\bet)  \ee
for all $r,s \in \ZZ$ and all $\al,\bet \in \Ga$. Regard $V$ as a
$\Ga'$--graded vector space (denoted by $V'$), such that, for all
$(r,\ga) \in \Ga'$,
 \bea  V'_{(1,\ga)} \eqq V_{\ga}  \\
       V'_{(r,\ga)} \eqq \Ze \;\;\;\mbox{if}\;\;\; r \neq 1 \,.  \eea
On the other hand, regard $L$ as a $\Ga'$--graded $\ve'$ Lie algebra (denoted
by $L'$), such that, for all $(r,\ga) \in \Ga'$,
 \bea  L'_{(0,\ga)} \eqq L_{\ga}  \\
       L'_{(r,\ga)} \eqq \Ze \;\;\;\mbox{if}\;\;\; r \neq 0 \,,  \eea
while leaving the product map unchanged. Obviously, $V'$ and $T(V)$ have a
natural structure of a $\Ga'$--graded $L'$--module.

A priori, $T(V')$ is a $\ZZ \ti \Ga' = \ZZ \ti \ZZ \ti \Ga$--graded vector
space, where the $\ZZ$--gradation corresponding to the first factor is the
canonical $\ZZ$--gradation of the tensor algebra $T(V')$, and where the one
corresponding to the second factor stems from the $\ZZ$--gradation of $V'$.
Obviously, these two $\ZZ$--gradations coincide. Thus we may simply forget
about one of them and regard $T(V')$ as a $\Ga'$--graded algebra and a
$\Ga'$--graded $L'$--module. A similar remark applies to the $\ve'$--symmetric
algebra $S(V',\ve')$ of $V'$.

With these conventions, $T(V)$ and $T(V')$ coincide, and $\Wept V$ is equal to
$S(V',\ve')$ (where these statements may be interpreted both in the sense of
$\Ga'$--graded algebras and of $\Ga'$--graded $L'$--modules).
 
The foregoing observation implies that for any result about $\ve'$--symmetric
algebras there is an analogous one for $\ve$--exterior algebras. In particular,
from Sec.~12 of Ref.~\cite{S.ma} we draw the following conclusions.

Let $n \ge 1$ be an integer and let
 \be  \om_n : V^n \lra \Wepfn V  \ee
be the canonical mapping defined by
 \be  \om_n(x_1,\ldots,x_n) = x_1 \wepf\ldots\wepf x_n  \ee
for all $x_i \in V$. Then $\om_n$ is an $\ve$--skew--symmetric $n$--linear
map, which is homogeneous of degree zero and $L$--invariant, and the pair
$(\Weptn V,\om_n)$ has the following universal property: \\
For any vector space $W$ and any $\ve$--skew--symmetric $n$--linear map
 \be  g : V^n \lra W \,,  \ee
there exists a unique linear map
 \be  \hat{g} : \Wepfn V \lra W  \ee
such that
 \be  g = \hat{g} \comp \om_n \,.  \ee
Now suppose in addition that $W$ is a graded $L$--module, and let
$\Lg^{a}_n(V,W;\ve)$ denote the graded $L$--submodule of $\Lgn(V,\ldots,V;W)$
consisting of all elements which are $\ve$--skew--symmetric. Then the
assignment $g \rar \hat{g}$ described above induces a graded $L$--module
isomorphism of $\Lg^{a}_n(V,W;\ve)$ onto $\Lg(\Weptn V,W)$.

Next, let $T^{a}_n(V,\ve)$ denote the graded $L$--submodule of $T_n(V)$
consisting of all $\ve$--skew--symmetric tensors in $T_n(V)$. The
$\ve$--skew--symmetrization
 \be  t \lra \frac{1}{n!}\sum_{\pi \in \Sn}\!\mbox{sgn}(\pi)\,S_{\pi}(t) \ee
defines a graded $L$--module homomorphism of $T_n(V)$ onto $T^{a}_n(V,\ve)$
(recall that the mappings $S_{\pi}$ have been defined in Ref.~\cite{S.ma}).
Using the universal property of $(\Weptn V,\om_n)$, it follows that there
exists a unique surjective graded $L$--module homomorphism
 \be  \chi_n : \Wepfn V \lra T^{a}_n(V,\ve)  \ee
such that
 \be  \chi_n(x_1 \wepf \ldots \wepf x_n) = 
\frac{1}{n!}\sum_{\pi\in\Sn}\!\mbox{sgn}(\pi)\,S_{\pi}(x_1 \ot\ldots\ot x_n) \ee
for all $x_i \in V$, and it is easy to see that $\chi_n$ is injective, and
hence an isomorphism of graded $L$--modules.

On the other hand, from the Poincar\'e, Birkhoff, Witt theorem \cite{S.gl}
we deduce that to any homogeneous basis of $V$ there corresponds, in the
obvious way, a homogeneous basis of $\Wept V$.

Finally, if $V$ and $W$ are two graded $L$--modules, then $\Wept(V \op W)$
(regarded as a $\ZZ \ti \Ga$--graded algebra and as a $\Ga$--graded
$L$--module) is canonically isomorphic to $(\Wept V) \otbeppr (\Wept W)$,
where $\otbeppr$ denotes the graded tensor product of $\ZZ \ti \Ga$--graded
algebras with respect to $\ve'$. This is a special case of an analogous
result for $\ve$ Lie algebras: If $G$ and $L$ are two $\ve$ Lie algebras, the
(universal) enveloping algebra $U(G \ti L)$ of the direct product $G \ti L$
of $G$ and $L$ is canonically isomorphic to the graded tensor product
$U(G) \otbep U(L)$ of the $\Ga$--graded algebras $U(G)$ and $U(L)$ (with
respect to $\ve$). The proof follows directly from the universal property of
enveloping algebras.

In particular, let $U = U_{\zb} \op U_{\ob}$ be a $\ZZ_2$--graded vector space.
Applying the preceding result to $V = U_{\zb}$ and $W = U_{\ob}$
(regarded as $\ZZ_2$--graded vector spaces), we conclude that the
super--exterior algebra $\Wept U$ is canonically isomorphic to
$(\bigwedge U_{\zb}) \bar{\ot} S(U_{\ob})$, where $\bigwedge U_{\zb}$ (resp.
$S(U_{\ob})$) is the {\em usual} exterior (resp. symmetric) algebra of the
vector space $U_{\zb}$ (resp. $U_{\ob}$), and where $\bar{\ot}$ denotes the
graded tensor product of $\ZZ$--graded algebras with respect to the commutation
factor $(r,s) \rar (-1)^{rs}$. In particular, we have the vector space
isomorphism
 \be  \Wepfn U \,\simeq\, \bop_{r=0}^{n}
                       ((\bigwedge^r U_{\zb}) \ot S_{n-r}(U_{\ob})) \,.  \ee

This enables us to make contact with the definition of the cohomology of Lie
superalgebras given in Ref.~\cite{Fu}. First, we notice that we have chosen a
different overall sign for the coboundary operator (which does not
matter). Moreover, there is another discrepancy, which (in our language)
amounts to dropping $\ga$ in the first sum on the right hand side of
Eq.~\reff{delex}. This is serious, for it spoils the $L$--invariance of
$\de^n$ (but not the $\ve$--skew-symmetry of $\de^n g$ and not the
validity of the fundamental equation \reff{deldel}). We suspect that
this is a misprint.

\sect{The tensorial square of the adjoint \newline
         $sl(1|2)$--module}
We keep the notation of Sec.~V and write $L$ instead of $sl(1|2)$. Our goal is
to investigate the graded $L$--module $L \ot L$\,. Of course, $L \ot L$
decomposes into the direct sum of the graded $L$--submodule $T^{s}_2(L,\ve)$
of super--symmetric tensors with the graded $L$--submodule $T^{a}_2(L,\ve)$
of super--skew--symmetric tensors,
 \be  L \ot L = T^{s}_2(L,\ve) \op T^{a}_2(L,\ve) \,.  \ee
Moreover, it is known that $T^{s}_2(L,\ve)$ (resp.~$T^{a}_2(L,\ve)$) is
isomorphic to $S_2(L,\ve)$ (resp.~to $\Weptt L$) (see Ref.~\cite{S.ma} and
App.~A).

The decomposition of $T^{a}_2(L,\ve)$ into the direct sum of simple graded
submodules has already been used in Sec.~V, we have
 \be  T^{a}_2(L,\ve) \simeq V(\ohal,\thhal,\ob) \op V(-\ohal,\thhal,\ob)
                                                    \op V(0,1,\zb) \,.  \ee
Thus all we have to do is to investigate $T^{s}_2(L,\ve)$. This $L$--module
is not semi--simple. In fact, we find that
 \be  T^{s}_2(L,\ve) \simeq V(0,2,\zb) \op V(0,1,\zb) \op V_8 \,, 
                                                             \label{dec} \ee
where $V_8$ is an 8--dimensional indecomposable graded $L$--module, which can
be described as follows. It has a basis consisting of four even weight vectors
$t,s,v,w$ and four odd weight vectors $\vpm , \wpm$\,, corresponding to the
weights $(0,0)$, $(0,0)$, $(1,0)$, $(-1,0)$ and $(\ohal,\pm\ohal)$,
$(-\ohal,\pm\ohal)$, respectively. The action of $L$ on $V_8$ is fixed by the
equations
 \bea   \Vpm \cd t = \vpm  \hst{1.5em}&,&\hst{2em}       \Wpm \cd t = \wpm \\
\Vpm \cd v_{\mp} = \pm v \hst{1.28em}&,&\hst{1.06em} \Wpm \cd w_{\mp} = \pm w \\
        \Wpm \cd v = \vpm  \hst{1.5em}&,&\hst{2em}       \Vpm \cd w = \wpm \eea
\vspace{-4.1ex}  
 \be  \Wpm \cd v_{\mp} = \Vpm \cd w_{\mp} = \pm s  \ee
 \be           \Vpm \cd s = \Wpm \cd s = 0 \,,  \ee
while the remaining images of the basis vectors under the action of $\Vpm$
or $\Wpm$ are equal to zero. We note that these equations fix the normalization
of the basis vectors up to an overall factor. The even vectors $t,s,v,w$ are
$sl(2)$--singlets, while $(\vpm)$ and $(\wpm)$ are standard $sl(2)$--doublets.
Obviously, $t$ generates the $L$--module $V_8$ (but, of course, it is not a
highest weight vector). Needless to say, $t$ can be replaced by any vector of
the form $t + \lam s$\,, with $\lam \in \KK$\,.

We do not want to specify the eight basis vectors inside $T^{s}_2(L,\ve)$.
Suffice it to note that $t$ can be chosen such that
 \be  t = Q_{+} \ot Q_{-} + Q_{-} \ot Q_{+} + 2\,Q_3 \ot Q_3
                                                     + 2\,B \ot B \,.  \ee
{\em Remark B.1.} The decomposition \reff{dec} of $T^s_2(L,\ve)$ has been
found independently of Ref.~\cite{Be}. As already noted, in this article
a decomposition of $U(L)$ into indecomposable $L$--submodules has been
determined, and \reff{dec} is just the first non--trivial part of it. The
module $V_8$ has been known since the early days of the theory of Lie
superalgebras. In fact, it was proved in Sec.~3.B of Ref.~\cite{SNR} that, for
both $\ga \in \ZZ_2$\,, we have
 \be  V(0,\ohal,\ga) \ot V(0,\ohal,\ga) \simeq V(0,1,\zb) \op V_8 \,, \ee
thus showing that the tensor product of two (typical) simple graded
$L$--modules need not be semi--simple.

The graded $L$--module $V_8$ has the following non--trivial graded
submodules:
 \bea       V_1 \eqq \KK\, s  \\
            V_4 \eqq \KK\, v_{+} \op \KK\, v_{-} \op \KK\, v \op \KK\, s  \\
           \Vfb \eqq \KK\, w_{+} \op \KK\, w_{-} \op \KK\, w \op \KK\, s  \\
            V_7 \eqq V_4 + \Vfb \,.  \eea
Obviously, $V_1$ is the submodule consisting of the $L$--invariant elements
of $V_8$\,. The modules $V_4$ and $\Vfb$ are indecomposable, and we have
 \be  V_{7}/\Vfb \simeq V_{4}/V_1 \simeq V(\ohal,\ohal,\ob)  \ee
 \be  V_{7}/V_4  \simeq \Vfb/V_1  \simeq V(-\ohal,\ohal,\ob) \,.  \ee
Finally, $V_{8}/V_7$ is a trivial one--dimensional $L$--module.

Let us next calculate $\mbox{dim}\,H^1(L,V_8)$\,. Using the information
obtained in Sec.~V, this is easily done by means of the long cohomology
sequence \reff{cohomseq}, applied to various short exact sequences of graded
$L$--modules. In fact, considering the sequence
 \be  \Ze \lra V_1 \lra V_4 \lra V_{4}/V_1 \lra \Ze  \ee
and the analogous one with $V_4$ replaced by $\Vfb$\,, we obtain
 \be  \mbox{dim}\,H^1(L,V_4) = \mbox{dim}\,H^1(L,\Vfb) = 1  \ee
 \be      H^2(L,V_4) = H^2(L,\Vfb) = \Ze \,.  \ee
Considering next the sequence
 \be  \Ze \lra V_4 \lra V_7 \lra V_{7}/V_4 \lra \Ze \,,  \ee
we conclude that
 \be  \mbox{dim}\,H^1(L,V_7) = 2 \,.  \ee
Finally, from the cohomology sequence corresponding to
 \be  \Ze \lra V_7 \lra V_8 \lra V_{8}/V_7 \lra \Ze  \ee
and from Eq.~\reff{hzero} we deduce that
 \be  \mbox{dim}\,H^1(L,V_8) = 1 \,. \label{hve}  \ee
The corresponding 1--cocycles can easily be determined. In fact, define two
linear maps
 \be  g : L \lra V_4 \;\;\;,\;\;\; \bar{g} : L \lra \Vfb  \ee
by
 \be  g(\Vpm) = \vpm \;\;\;,\;\;\; \bar{g}(\Wpm) = \wpm  \ee
 \be  -g(B) = \bar{g}(B) = s \,,  \ee
with $g(X) = \bar{g}(Y) = 0$ for all other elements $X,Y$ of our standard basis
of $L$\,. Then $g$ and $\bar{g}$ are 1--cocycles with values in $V_4$ and
$\Vfb$\,, respectively, and they are not coboundaries. If we regard $g$ and
$\bar{g}$ as 1--cocycles with values in $V_7$\,, their cohomology classes in
$H^1(L,V_7)$ are linearly independent. Finally, if we regard $g$ and $\bar{g}$
as 1--cocycles with values in $V_8$\,, they are still not coboundaries, but
$g$ and $-\bar{g}$ are cohomologous, with
 \be  g + \bar{g} = \de^{0}t \,.  \ee

The foregoing results have an immediate bearing on the cohomology of $U(L)$.
In fact, it is known \cite{S.ce} that $U(L)$, regarded as a graded $L$--module
(and also as a graded coalgebra) is canonically isomorphic to $S(L,\ve)$.
Consequently, there exists a graded $L$--submodule of $U(L)$, which is a
direct summand of $U(L)$ and which is isomorphic to $V_8$\,. (Note that the
simple graded subquotients of $V_8$ are exactly those which according to
Lemma 5.1 are allowed.) Using the isomorphism \reff{dirsum}, we conclude
from Eq.~\reff{hve} that
 \be  H^1(L,U(L)) \neq \Ze \,.  \ee

\end{appendix}

\end{document}